\documentclass[a4paper,12pt]{amsart}
\usepackage[T2A]{fontenc}
\usepackage[cp1251]{inputenc}
\usepackage[english]{babel}
\usepackage{amssymb, amsmath, latexsym,mathrsfs}
\usepackage{amsfonts,mathtext,cite,enumerate,float,amsthm}
\usepackage{color}
\numberwithin{equation}{section}
\usepackage{graphicx}
\renewcommand{\le}{\leqslant}
\renewcommand{\ge}{\geqslant}
\newcommand{\diam}{\text{\,\rm diam\,}}
\newcommand{\dist}{\text{\,\rm dist\,}}
\renewcommand{\Re}{\text{\rm Re\,}}
\renewcommand{\Im}{\text{\rm Im\,}}

\newtheoremstyle{citedth}%
  {5pt}
  {5pt}
  {\itshape}
  {}
  {\bfseries}
  {.}
  {.3em}
  {\thmname{#1} \thmnumber{#2} \thmnote{\normalfont#3}}

\theoremstyle{citedth}
\newtheorem{theoremA}{Theorem}

\theoremstyle{theorem}
\newenvironment{Proof}
{\noindent{\bf Proof.}}
{\hfill$\scriptstyle\blacksquare$}
\newtheorem{theorem}{Theorem}
\newtheorem{proposition}{Proposition}
\newtheorem*{theorem*}{Theorem~\ref{new_theorem2}}
\newtheorem{lemma}{Lemma}
\theoremstyle{definition}
\newtheorem{example}{Example}
\newtheorem{remark}{Remark}

\usepackage{geometry}
\begin{document}
\title[Operators whose $L^2$-boundedness implies rectifiability]{A new family of singular integral operators whose $L^2$-boundedness implies rectifiability}
\author{Petr Chunaev}
\thanks{The research was supported by the ERC grant 320501 of the European Research Council (FP7/2007-2013).}
\address{Departament de Matem\`{а}tiques, Universitat Aut\`{o}noma de Barcelona, Edifici C
08193 Bellaterra (Barcelona) Spain} \email{chunayev@mail.ru,
chunaev@mat.uab.cat}

\begin{abstract}
Let $E \subset \mathbb{C}$ be a Borel set such that
$0<\mathcal{H}^1(E)<\infty$. David and L\'{e}ger proved that the
Cauchy kernel $1/z$ (and even its coordinate parts $\Re z/|z|^2$ and
$\Im z/|z|^2$, $z\in \mathbb{C}\setminus\{0\}$) has the following property: the $L^2(\mathcal{H}^1\lfloor E)$-boundedness of the
corresponding singular integral operator implies that $E$ is rectifiable. Recently Chousionis, Mateu, Prat and Tolsa extended this
result to any kernel of the form $(\Re z)^{2n-1}/|z|^{2n}$, $n\in
\mathbb{N}$. In this paper, we prove that the above-mentioned property holds for operators associated with the much wider class of the kernels
$(\Re z)^{2N-1}/|z|^{2N}+t\cdot(\Re z)^{2n-1}/|z|^{2n}$, where
$n$ and $N$ are positive integer numbers such that $N\ge n$, and $t\in \mathbb{R}\setminus (t_1,t_2)$ with $t_1,t_2$ depending only on $n$ and~$N$.
\end{abstract}

\maketitle

\section{Introduction}
Let $E\subset \mathbb{C}$ be a Borel set and $B(z,r)$ an open disc with center $z\in \mathbb{C}$ and radius~$r>0$. By $\mathcal{H}^1(E)$ we denote the one-dimensional Hausdorff measure of $E$. A~set $E$ is called \textit{rectifiable} if it is contained, up to an $\mathcal{H}^1$-negligible set, in a countable union of Lipschitz graphs. A set $E$ with $\mathcal{H}^1(E)<\infty$ is called \textit{purely unrectifiable} if it intersects any Lipschitz graph in a set of $\mathcal{H}^1$-measure zero.

A singular integral operator $T_K$ associated with a kernel $K$ is
formally defined as
\begin{equation}
\label{integral_operator}
T_Kf(z):=\int_E f(\zeta)K(z-\zeta)\,d\mathcal{H}^1(\zeta),
\end{equation}
where $K:\mathbb{C}\setminus\{0\}\to \mathbb{C}$ is a standard
kernel (see its definition, for instance, in \cite{Christ}) and $f$
is some reasonable function, say, $f\in L^1(\mathcal{H}^1\lfloor
E)$. The integral in (\ref{integral_operator}) might not converge
absolutely and therefore one usually considers $T_{K,\varepsilon}$, a truncated version
of $T_K$, which is defined by the
above-mentioned integral but over the set $E \setminus
B(z,\varepsilon)$ for some $\varepsilon>0$. The operator $T_K$ is
said to be \textit{bounded on $L^2(\mathcal{H}^1\lfloor E)$} if the
operators $T_{K,\varepsilon}$ are bounded on
$L^2(\mathcal{H}^1\lfloor E)$ uniformly on $\varepsilon$. We also
recall that \textit{the~principal value} (\textrm{p.v.}) of the
operator $T_K$ is said to exist $\mathcal{H}^1$-a.e. on $E$ if
$\lim_{\varepsilon\to 0+}T_{K,\varepsilon}f(z)$ exists and is finite
for almost every $z\in E$ and $f$ from a reasonable functional
space.

The connection between the $L^2$-boundedness, existence of p.v. of
$T_K$ (defined for more general measures than $\mathcal{H}^1$) on a
set $E$ and the geometric properties of this set, e.g.
rectifiability, is an object of intensive investigations. They were
initiated by Calder\'{o}n \cite{Calderon}, who proved that the
Cauchy transform, i.e. $T_K$ with $K(z)=1/z$, is $L^2$-bounded on
Lipschitz graphs with small slope. Later on, Coifman, McIntosh and
Meyer \cite{CMM} removed the small Lipschitz constant assumption. In
\cite{David} David fully characterized rectifiable \textit{curves}
$\Gamma$, for which the Cauchy transform is bounded on
$L^2(\mathcal{H}^1\lfloor \Gamma)$: they have to satisfy the linear
growth condition
$$
\mathcal{H}^1(\Gamma\cap B(z,r))\le Cr,\qquad r>0,\qquad z\in \mathbb{C}.
$$
These results led to further development of tools
for understanding the above-mentioned connection.
 For more information about this topic, see the corresponding parts of \cite{Mattila2002,Mattila1997,Tolsa_book, David_Semmes,Mattila}.

A second wind in the area happened after the discovery of the so-called curvature method, which became very influential in the study of the Cauchy transform and analytic capacity \cite{M,David1998,Tolsa2003,MMV}. We now describe the core of the method. Given pairwise distinct points $z_1,z_2,z_3\in \mathbb{C}$, their \textit{Menger curvature} is
$$
c(z_1,z_2,z_3)=\frac{1}{R(z_1,z_2,z_3)},
$$
where $R(z_1,z_2,z_3)$ is the radius of the circle passing through
$z_1$, $z_2$ and $z_3$
 (with $R(z_1,z_2,z_3)=\infty$ and $c(z_1,z_2,z_3)=0$ if the points are collinear). It is easily seen that the curvature can
  be calculated in geometrical terms in different ways, e.g.
\begin{equation}
\label{curvature_geometry}
c(z_1,z_2,z_3)=\frac{4S(z_1,z_2,z_3)}{|z_1-z_2||z_1-z_3||z_2-z_3|}=\frac{2\sin \widehat{z_1z_2z_3}}{|z_1-z_3|},
\end{equation}
where $S(z_1,z_2,z_3)$ stands for the area of the triangle
$(z_1,z_2,z_3)$ and $\widehat{z_1z_2z_3}$ is the angle of this
triangle opposite to the side $z_1z_3$.

The relationship between the curvature and the Cauchy kernel
originates from the following identity due to Melnikov \cite{M}:
\begin{equation}
\label{c^2}
 c(z_1,z_2,z_3)^2=\sum_{s\in
\mathfrak{S}_3}\frac{1}{(z_{s_2}-z_{s_1})\overline{(z_{s_3}-z_{s_1})}},
\end{equation}
where $\mathfrak{S}_3$ is the group of permutations of three
elements. It is very important that the quantity in the right hand side turns out to be a non-negative real number.  Let us also define the so-called \textit{curvature of a Borel measure} $\mu$:
\begin{equation}
\label{c2(mu)}
 c^2(\mu)=\iiint
c(z_1,z_2,z_3)^2\;d\mu(z_1)\,d\mu(z_2)\,d\mu(z_3),
\end{equation}
introduced by Melnikov, too.
One can consider $c^2_{\varepsilon}(\mu)$, a truncated version of $c^2(\mu)$, which is the above-mentioned triple integral over the set
$$
\{(z_1,z_2,z_3)\in \mathbb{C}^3: |z_k-z_j|\ge \varepsilon>0, \quad 1\le k,j\le 3, \quad j\neq k\}.
$$
Clearly, $c^2_{\varepsilon}(\mu)$ is always
non-negative due to~(\ref{c^2}).

If $\mu$ is a finite Borel measure with linear growth, i.e.
$\mu(B(z,r))\le Cr$ for all $z\in \textrm{spt}\,\mu$, then the
relation between the curvature and the $L^2(\mu)$-norm of the Cauchy
transform (of $f(\zeta)\equiv 1$) is specified by the following Melnikov-Verdera
identity \cite{MV}:
\begin{equation}
\label{Melnikov-Verdera}
\int\left|\int_{\mathbb{C}\setminus B(z,\varepsilon)}\frac{d\mu(\zeta)}{\zeta-z}\right|^2d\mu(z)=
\tfrac{1}{6}c_\varepsilon^2(\mu)+O(\mu(\mathbb{C})),\qquad |O(\mu(\mathbb{C}))|\le C\mu(\mathbb{C}).
\end{equation}

Relying on the curvature method described above, Mattila, Melnikov and Verdera \cite{MMV} proved that if $E \subset \mathbb{C}$ is an \textit{Ahlfors-David regular
set}, i.e. ${0<\mathcal{H}^1(E)<\infty}$ and
there exists a constant $C>0$ such that
$$
C^{-1}r\le \mathcal{H}^1(E\cap B(z,r))\le Cr, \qquad z\in E, \qquad 0<r< \mathcal{H}^1(E),
$$
then the Cauchy transform is $L^2(\mathcal{H}^1\lfloor E)$-bounded
if and only if $E$ is contained in an Ahlfors-David regular curve.

Later on, essentially using that the curvature is non-negative, David and L\'{e}ger made the
following deep contribution (see a brief exposition of the proof in
Section~\ref{section_proof}).
\begin{theoremA}[\cite{L}]
\label{theorem_Leger}
Given a Borel set $E\subset\mathbb{C}$ such that
$0<\mathcal{H}^1(E)<\infty$, if ${c^2(\mathcal{H}^1\lfloor
E)<\infty}$, then $E$ is rectifiable. Moreover, if the Cauchy
transform is $L^2(\mathcal{H}^1\lfloor E)$-bounded, then $E$ is
rectifiable.
\end{theoremA}

This theorem is stated in \cite{L} for $\mathbb{R}^d$ but we formulate it only for $\mathbb{R}^2=\mathbb{C}$.
Note also that its second statement is a direct corollary of the first
one and (\ref{Melnikov-Verdera}).

Until recently, very few things were known in this direction beyond
the Cauchy kernel. For instance, the same result is true for its
coordinate parts \cite{CMPT}:
\begin{equation}
\label{coordinate_parts}
\Re z/|z|^2 \quad \text{or}\quad \Im z/|z|^2,\qquad z\in \mathbb{C}\setminus\{0\}.
\end{equation}
Indeed, consider the following permutations:
\begin{equation}
\label{permutation_Main}
\begin{split}
p_K(z_1,z_2,z_3):=&K(z_{1}-z_{2})K(z_{1}-z_{3})+K(z_{2}-z_{1})K(z_{2}-z_{3})\\&+K(z_{3}-z_{1})K(z_{3}-z_{2}),
\end{split}
\end{equation}
where $K$ is some standard real kernel. Then it is not difficult to show
that if $K$ is one of the kernels (\ref{coordinate_parts}), then
\begin{equation}
\label{permutation_real_part}
p_{K}(z_1,z_2,z_3)=\tfrac{1}{4}c(z_1,z_2,z_3)^2.
\end{equation}
This fact was a motivation point of the recent paper \cite{CMPT} by Chousionis, Mateu, Prat and Tolsa. The David-L\'{e}ger result is extended there to the kernels
\begin{equation}
\label{kappa_n}
\kappa_n(z):=\frac{(\Re z)^{2n-1}}{|z|^{2n}}, \qquad n\in \mathbb{N}.
\end{equation}
Namely, it is shown in \cite{CMPT} that for any given triple
$(z_1,z_2,z_3)\in \mathbb{C}^3$,
\begin{equation}
\label{positive_permut_kappa_m} p_{\kappa_n}(z_1,z_2,z_3)\ge 0,
\end{equation}
and $p_{\kappa_n}(z_1,z_2,z_3)=0$ if and only if the points
$z_1,z_2,z_3$ are collinear. Moreover, it is proved that the
permutations $p_{\kappa_n}(z_1,z_2,z_3)$ behave similarly to  $c^2(z_1,z_2,z_3)$ for triangles with comparable
sides, whose one side makes a big angle with the vertical line. This fact
enables the authors of \cite{CMPT} to adapt the method from \cite{L}
to the kernels $\kappa_n$. This adaptation however requires them to make
several essential modifications in crucial points, where the curvature must be exchanged for the
permutations $p_{\kappa_n}$, and provide new arguments whenever the scheme of L\'{e}ger does not work
(see also Section~\ref{section_proof} for more details).

To state the corresponding result we need the following generalization of~(\ref{c2(mu)}):
\begin{equation}
\label{p_K(mu)}
p_K(\mu)=\iiint p_K(z_1,z_2,z_3)\;d\mu(z_1)\,d\mu(z_2)\,d\mu(z_3).
\end{equation}
\begin{theoremA}[\cite{CMPT}]
\label{CMPT_theorem}
Let $n\ge 1$. Given a Borel set $E\subset
\mathbb{C}$ such that $0<\mathcal{H}^1(E)<\infty$,
if~$p_{\kappa_n}(\mathcal{H}^1\lfloor E)<\infty$, then $E$ is
rectifiable. Moreover, if the operator $T_{\kappa_n}$ is
$L^2(\mathcal{H}^1\lfloor E)$-bounded, then $E$ is rectifiable.
\end{theoremA}

Obviously, for $n=1$ one gets the statement for the real part of the
Cauchy kernel. However, for $n\ge 2$
this is the first example of singular integral
operators with the above-mentioned property, which are not directly related
to the Cauchy transform.

In this paper, we consider a linear combination of the kernels
(\ref{kappa_n}), namely, the parametric kernels
\begin{equation}
\label{kernels_n_N} K_t(z):=\kappa_N(z)+t\cdot \kappa_n(z), \qquad
N>n, \qquad n,N\in \mathbb{N^+},\qquad t\in \mathbb{R},
\end{equation}
where the parameters $n$, $N$ and $t$ are fixed. We find
values of $t$, depending on $n$ and $N$, such that a result,
analogous to the David-L\'{e}ger theorem, is valid for singular
integral operators associated with the kernels $K_t$. For this
purpose we first study the sign of the permutations
(\ref{permutation_Main}) for the kernels (\ref{kernels_n_N}) and then,
for the case when these permutations are non-negative, adapt the scheme from \cite{CMPT} to prove the result of David-L\'{e}ger type. The next section contains the corresponding statements.

\section{Main results}

\begin{figure}
  \includegraphics[width=15cm]{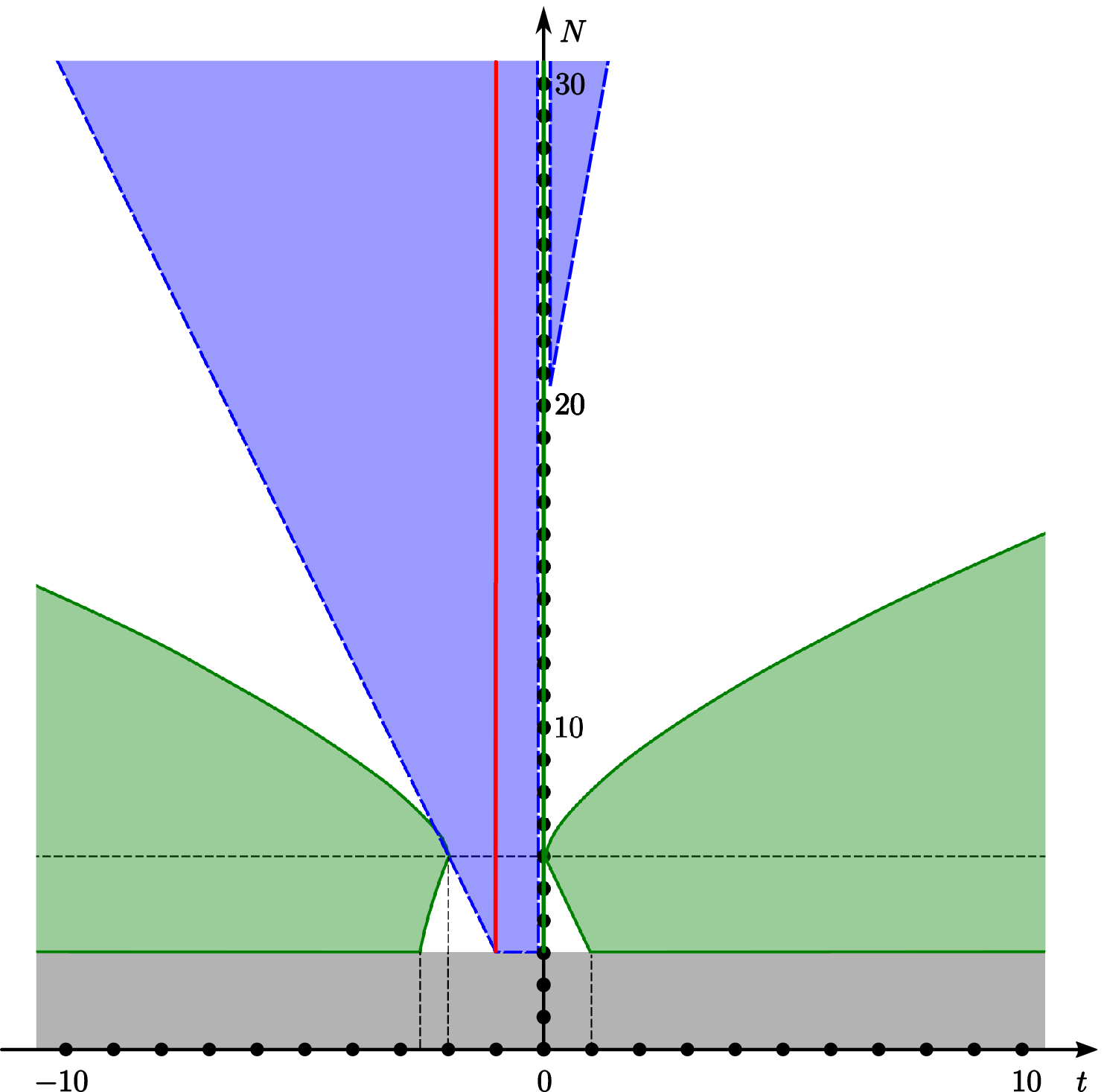}\\
 \caption{}
  \label{FIG1}
\end{figure}

First of all let us mention that the case $t=0$ in the theorems below agrees with the inequality
(\ref{positive_permut_kappa_m}) and Theorem~\ref{CMPT_theorem},
proved in \cite{CMPT}. We now indicate the values of $t$ such that the
permutations $p_{K_t}(z_1,z_2,z_3)$ are non-negative for all triples
$(z_1,z_2,z_3)$.
\begin{theorem}
\label{new_theorem1} Let $K_t$ be a kernel of the form
$(\ref{kernels_n_N})$ with $t=0$ or
\begin{eqnarray}
\label{estimates_t_N<2n}
&t\in \mathbb{R}\setminus \left(-\tfrac{1}{2}\left(3+\sqrt{9-4\tfrac{N}{n}}\right); 2-\tfrac{N}{n}\right),\qquad & n<N\le 2n,\\
\label{estimates_t_N>2n}
&t\in \mathbb{R}\setminus \left(-\frac{1}{2}\left(3+\rho_{n,N}+\sqrt{(3+\rho_{n,N})^2-4\frac{N}{n}}\right); \rho_{n,N}\right),\quad &N\ge 2n,
\end{eqnarray}
where $\rho_{n,N}:=(\tfrac{N}{n}-2)\sqrt{N-2n}$. In particular, $t\in \mathbb{R}\setminus (-2;0)$ for $N=2n$.
Then
$$
p_{K_t}(z_1,z_2,z_3)\ge 0 \quad \text{for all} \quad
(z_1,z_2,z_3)\in \mathbb{C}^3.
$$

Furthermore, the range of the parameter $t$ in the case $N=2n$ is sharp.
\end{theorem}
\begin{remark}
\label{remark1} The conditions (\ref{estimates_t_N<2n}) and
(\ref{estimates_t_N>2n}), guaranteeing that
$p_{K_t}(z_1,z_2,z_3)\ge 0$, cannot be weakened much in the following
sense. As we will show in Section~\ref{section_examples}, there are
triples $(z_1,z_2,z_3)$ such that $p_{K_t}(z_1,z_2,z_3)$ can change
 sign if
\begin{equation}
\label{t_bounds}
\begin{array}{ll}
  t\in(-\frac{N}{n};0)  & \text{for all }n \text{ and } N,\phantom{\frac{L}{Q_{q_q}}} \\
  t\in\left(0;\frac{2}{e^{3/2}}\frac{N}{n}\right)  & \text{for }N\gg n \quad (\text{i.e. } N \text{ is large enough with respect to } n).
\end{array}
\end{equation}
Surprisingly, in this context the case $t=0$ is an isolated point
with non-negative permutations. Thus the curvature method,
requiring the permutations to be non-negative, cannot be applied
directly for $t$ indicated in (\ref{t_bounds}).

From the aforesaid it follows that (\ref{estimates_t_N<2n}) and
(\ref{estimates_t_N>2n}) are indeed sharp for $N=2n$.

Figure~\ref{FIG1} illustrates Theorem~\ref{new_theorem1} and this remark for $n=3$ and different $N>3$
(consider the horizontal line passing through a fixed
\textit{positive integer} $N$ in order to determine the
corresponding $t$). The green area represents the values of $t$,
indicated in (\ref{estimates_t_N<2n}) and (\ref{estimates_t_N>2n}),
i.e. those where  $p_{K_t}$ are guaranteed to be non-negative (the
boundaries are included). In the blue area (the boundaries are not
included), obtained by computer-based exhaustive search, the
permutations can change sign. Note that the part of the blue area
for $t<0$ is exactly the former in (\ref{t_bounds}). Moreover, the
part for $t>0$ quite agrees with the latter in (\ref{t_bounds}). The
white area is not covered by our results and, generally speaking, we
can say nothing about the sign of $p_{K_t}$ therein. However,
computer experiments suggest that the permutations $p_{K_t}$ are
non-negative everywhere except the blue area and thus (\ref{t_bounds})
seems to give likely boundaries for~$t$, whose corresponding
permutations can change sign.
\end{remark}

Relying deeply on Theorem~\ref{new_theorem1}, we will prove the
following David-L\'{e}ger type result.
\begin{theorem}
\label{new_theorem2}
Let $K_t$ be a kernel of the form $(\ref{kernels_n_N})$ with $t$, mentioned in Theorem \ref{new_theorem1}.
Given a Borel set $E\subset \mathbb{C}$ such that $0<\mathcal{H}^1(E)<\infty$,
if $p_{K_t}(\mathcal{H}^1\lfloor E)<\infty$, then $E$ is rectifiable.
Moreover, if the operator $T_{K_t}$ is $L^2(\mathcal{H}^1\lfloor E)$-bounded, then $E$ is rectifiable.
\end{theorem}
\begin{remark}
\label{remark2}
 It is known that for $t=-1$ (see the red line in
Figure~\ref{FIG1}), which belongs to the area, where the
permutations can change sign,  the statement of
Theorem~\ref{new_theorem2} is not valid anymore, i.e.
$L^2$-boundedness  does not imply rectifiability. Indeed,  it is a
corollary of the following result due to Huovinen (note that all
kernels $K_t$ of the form (\ref{kernels_n_N}) with $t=-1$  belong to
Huovinen's class $\mathscr{H}$).
\begin{theoremA}[\cite{H}]
Let $K$ belong to the class $\mathscr{H}$ of odd kernels satisfying
\begin{eqnarray*}
|K(x-y)-K(x-z)|\le \frac{C\; |y-z|}{|x-y||x-z|}, \qquad |K(z)|\le \frac{C}{|z|},\qquad x,y,z\in \mathbb{C},\\
K(r)=0,\qquad K(z)=-K(-\overline{z}),\qquad r\in \mathbb{R},\qquad C={\rm const}.
\end{eqnarray*}
Then there exists a purely unrectifiable set $E$ with
$\mathcal{H}^1(E)>0$, such that the operator $T_K$, associated with
the kernel $K$,
 is bounded on $L^2(\mathcal{H}^1\lfloor E)$ and, moreover, ${\rm p.v.}\;T_{K}$ exists $\mathcal{H}^1$-a.e. and is finite.
\end{theoremA}

For a particular family of the kernels (\ref{kernels_n_N}) with
$(n,N)=(1,2)$, i.e.
$$
k_t(z):=\frac{(\Re z)^{3}}{|z|^{4}}+t\cdot\frac{\Re z}{|z|^{2}},
$$
even more is known. Namely,  for $t=-3/4$, whose corresponding $p_{K_t}$
can change sign, there also exists a purely
unrectifiable set $E$ such that $T_{k_t}$ is
$L^2(\mathcal{H}^1\lfloor E)$-bounded. One can get this from the
following result of Jaye and Nazarov for the kernel
$k(z):=\overline{z}/z^2$, noting that $\Re k(z)=4k_t(z)$ for
$t=-3/4$. We formulate it in a slightly different form than in
\cite{JN}.
\begin{figure}
  \includegraphics[width=15cm]{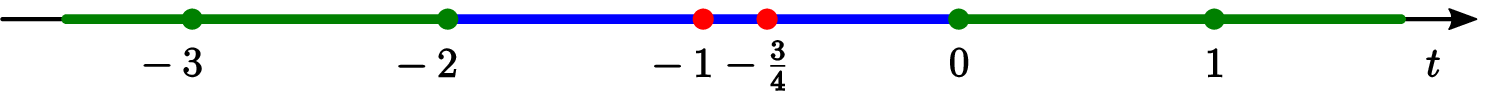}\\
 \caption{}
  \label{FIG2}
\end{figure}
\begin{theoremA}[\cite{JN}]
There exists a purely unrectifiable set $E$ with
$\mathcal{H}^1(E)>0$, such that the operator $T_k$, associated with
the kernel $k$, is bounded on $L^2(\mathcal{H}^1\lfloor E)$ and,
moreover,  ${\rm p.v.}\;T_k$ fails to exist $\mathcal{H}^1$-a.e.
\end{theoremA}
\end{remark}

Figure~\ref{FIG2} illustrates known results for the kernels
$k_t$. By Theorems~\ref{new_theorem1} and~\ref{new_theorem2}, if
$t\in \mathbb{R}\setminus (-2;0)$, then the
permutations $p_{k_t}$ are non-negative and the
$L^2(\mathcal{H}^1\lfloor E)$-boundedness of $T_{k_t}$ implies the
rectifiability of $E$ (see the intervals coloured green). By the arguments from Remarks~\ref{remark1}
and~\ref{remark2}, the permutations $p_{k_t}$ for $t\in (-2;0)$ change sign (the interval coloured blue) and there are two values of $t$ (the red points) such that the operator $T_{K_t}$ is
$L^2(\mathcal{H}^1\lfloor E)$-bounded but $E$ is not rectifiable.

\section{Auxiliary results}
\label{section-auxiliary}

This section is devoted to several auxiliary lemmas, which will be used to prove Theorems~\ref{new_theorem1}~and~\ref{new_theorem2} in Section~\ref{section_proof}.

If a kernel $K$ is real and odd, then one can show that the
permutations (\ref{permutation_Main}) are invariant under
translations. This can be done, e.g. by  the substitutions $u=z_1-z_2$ and $v=z_1-z_3$.
Consequently, one point can be always fixed and it is enough to
consider only permutations of the form
\begin{equation}
\label{p_K} p_{K}(0,u,v)=K(u)K(v)+K(u)K(u-v)+K(v)K(v-u),
\end{equation}
where $u,v\in \mathbb{C}\setminus \{0\}$ are distinct points.
The kernels (\ref{kappa_n}) and (\ref{kernels_n_N}) that we study are real and odd and hence we can use
(\ref{p_K}) instead of (\ref{permutation_Main}). Furthermore, the case of collinear points $u$ and $v$ is trivial as then $p_K(0,u,v)\equiv 0$ and thus we can skip it.

We will use the following lemma many times below. Note that it can be easily generalized for any other couple of kernels instead of $\kappa_n$ and $\kappa_N$.
\begin{lemma}
\label{lemma1}
Given $K_t$ of the form $(\ref{kernels_n_N})$,
\begin{equation}
\label{parabola_t}
p_{K_t}(0,u,v)=p_{\kappa_N}(0,u,v)+\varphi_{n,N}(0,u,v)\,t+p_{\kappa_n}(0,u,v)\,t^2,
\end{equation}
where
\begin{equation}
\label{varphi}
\begin{array}{ll}
\varphi_{n,N}(0,u,v)=&\kappa_N(u)(\kappa_n(v)+\kappa_n(u-v))\\
&  +\kappa_N(v)(\kappa_n(u)+\kappa_n(v-u))
                        +\kappa_N(u-v)(\kappa_n(u)-\kappa_n(v)).
\end{array}
\end{equation}
\end{lemma}
\begin{Proof}
We substitute (\ref{kernels_n_N}) into (\ref{p_K}) and get
$$
\begin{array}{l}
p_{K_t}(0,u,v)\\
=\kappa_N(u)\kappa_N(v)+(\kappa_N(u)\kappa_n(v)+\kappa_n(u)\kappa_N(v))\,t+\kappa_n(u)\kappa_n(v)\,t^2\\
\;
+\kappa_N(u)\kappa_N(u-v)+(\kappa_N(u)\kappa_n(u-v)+\kappa_n(u)\kappa_N(u-v))\,t+\kappa_n(u)\kappa_n(u-v)\,t^2\\
\;
+\kappa_N(v)\kappa_N(v-u)+(\kappa_N(v)\kappa_n(v-u)+\kappa_n(v)\kappa_N(v-u))\,t+\kappa_n(v)\kappa_n(v-u)\,t^2.
\end{array}
$$
To finish the proof it is enough to group the terms and take into account (\ref{p_K}).
\end{Proof}

It is important that the leading coefficient of the quadratic polynomial
(\ref{parabola_t}) (with respect to $t$) is always non-negative by the property (\ref{positive_permut_kappa_m}).

From now on, in order to simplify formulas we skip $(0,u,v)$ in
permutations and other expressions if there is no confusion. For example, we
write $p_K$ instead of $p_K(0,u,v)$. In addition, we use the following
notations:
\begin{equation}
\label{notation} \lambda_1:=\frac{\Re u}{|u|},\qquad
\lambda_2:=\frac{\Re v}{|v|},\qquad \lambda_3:=\frac{\Re
(u-v)}{|u-v|}, \qquad \Lambda:=\lambda_1\lambda_2\lambda_3,
\end{equation}
where the denominators do not vanish as the points $u$ and $v$ are assumed
to be distinct and non-collinear. Note that in these terms,
\begin{equation}
\label{p_repr}
p_{\kappa_n}=\frac{(\lambda_1\lambda_2)^{2n-1}}{|u||v|}+\frac{(\lambda_1\lambda_3)^{2n-1}}{|u||u-v|}
-\frac{(\lambda_2\lambda_3)^{2n-1}}{|v||v-u|}
\end{equation}
and
\begin{align}
\label{phi_repr} \varphi_{n,N}=&\frac{\lambda_1^{2N-1}}{|u|}\left(\frac{\lambda_2^{2n-1}}{|v|}+\frac{\lambda_3^{2n-1}}{|u-v|}\right)\\
\nonumber
&+\frac{\lambda_2^{2N-1}}{|v|}\left(\frac{\lambda_1^{2n-1}}{|u|}-\frac{\lambda_3^{2n-1}}{|v-u|}\right)
+\frac{\lambda_3^{2N-1}}{|u-v|}\left(\frac{\lambda_1^{2n-1}}{|u|}-\frac{\lambda_2^{2n-1}}{|v|}\right).
\end{align}

What is more, another representation of $\varphi_{n,N}$ is valid.
\begin{lemma}
\label{lemma2}
In terms of $(\ref{notation})$ it holds that
\begin{equation}
\label{phi_K1_K2}
\varphi_{n,N}=\tau_1 p_{\kappa_n}-\tau_2,
\end{equation}
where
\begin{equation}
\label{gamma1} \tau_1:=\lambda_1^{2(N-n)}+\lambda_2^{2(N-n)}+\lambda_3^{2(N-n)},\qquad 0\le \tau_1\le 3,
\end{equation}
and
\begin{equation}
\label{gamma2} \tau_2:=\Lambda^{2(N-n)}
\left(\frac{(\lambda_1\lambda_2)^{2(2n-N)-1}}{|u||v|}+\frac{(\lambda_1\lambda_3)^{2(2n-N)-1}}{|u||u-v|}
-\frac{(\lambda_2\lambda_3)^{2(2n-N)-1}}{|v||v-u|}\right).
\end{equation}
In particular, $\tau_2\equiv 0$ if $N=2n$.
\end{lemma}
\begin{Proof}
Direct multiplication of $\tau_1$ by $p_{\kappa_n}$ gives
$$
\left(\lambda_1^{2(N-n)}+\lambda_2^{2(N-n)}+\lambda_3^{2(N-n)}\right)\cdot
\left(\frac{(\lambda_1\lambda_2)^{2n-1}}{|u||v|}+\frac{(\lambda_1\lambda_3)^{2n-1}}{|u||u-v|}
-\frac{(\lambda_2\lambda_3)^{2n-1}}{|v||v-u|}\right)
$$
$$
=\left(
\frac{\lambda_3^{2(N-n)}(\lambda_1\lambda_2)^{2n-1}}{|u||v|}
+\frac{\lambda_2^{2(N-n)}(\lambda_1\lambda_3)^{2n-1}}{|u||u-v|}
-\frac{\lambda_1^{2(N-n)}(\lambda_2\lambda_3)^{2n-1}}{|v||v-u|}
\right)
$$
$$
+\frac{\lambda_1^{2N-1}}{|u|}\left(\frac{\lambda_2^{2n-1}}{|v|}+\frac{\lambda_3^{2n-1}}{|u-v|}\right)
+\frac{\lambda_2^{2N-1}}{|v|}\left(\frac{\lambda_1^{2n-1}}{|u|}-\frac{\lambda_3^{2n-1}}{|v-u|}\right)
+\frac{\lambda_3^{2N-1}}{|u-v|}\left(\frac{\lambda_1^{2n-1}}{|u|}-\frac{\lambda_2^{2n-1}}{|v|}\right),
$$
which is exactly $\tau_2+\varphi_{n,N}$ by (\ref{phi_repr}) and
(\ref{gamma2}).
\end{Proof}
\begin{lemma}
\label{lemma_fresh}
Given $\kappa_n$ and $\kappa_N$ of the form $(\ref{kappa_n})$,
\begin{equation}
\label{connection_kernels_n_and_m}
\frac{N}{n}\cdot \Lambda^{2(N-n)}\cdot p_{\kappa_n}\le p_{\kappa_N}, \qquad 1\le n \le N.
\end{equation}
\end{lemma}
Note that  this inequality for $n=1$ was obtained in \cite[Proof of Lemma~2.3]{CMPT}.
We will use the following lemma from there in order to prove the general form.
\begin{lemma}[Proof of Proposition 2.1 in \cite{CMPT}]
\label{lemma_CMPT}
One has the representation
$$
p_{\kappa_m}=\sum_{k=1}^m \binom{m}{k} \Lambda^{2(m-k)} h_k(u,v),
$$
where $h_k(u,v)\ge 0$ and are defined as follows:
\begin{align*}
h_k(u,v)=&(|u||v||u-v|)^{-2k}\left((\Re u\;\Re v)^{2k-1}(\Im (u-v))^{2k}\right.\\
&\left.+(\Re u\;\Re(u-v))^{2k-1}(\Im v)^{2k}+(\Re
v\;\Re(v-u))^{2k-1}(\Im u)^{2k}\right).
\end{align*}
\end{lemma}
\begin{Proof}
Within the settings of Lemma~\ref{lemma_CMPT},
$$
  \Lambda^{2(N-n)}\frac{p_{\kappa_n}}{p_{\kappa_N}}
 = \frac{\sum_{k=1}^{n} \binom{n}{k}\binom{N}{k}^{-1}H_{k}(u,v)}{\sum_{k=1}^N H_{k}(u,v)},
$$
where $H_{k}(u,v):=\binom{N}{k}\Lambda^{2(N-k)}h_k(u,v)\ge 0$. Furthermore,
$$
\binom{n}{k}\binom{N}{k}^{-1}=\frac{n!}{(n-k)!}\frac{(N-k)!}{N!}=\frac{(n-k+1)\cdots
n}{(N-k+1)\cdots N}\le \frac{n}{N},\qquad 1\le k\le n,
$$
and finally
\begin{equation*}
\label{interest_estimate}
\Lambda^{2(N-n)}\frac{p_{\kappa_n}}{p_{\kappa_N}}\le \frac{n}{N}\cdot \frac{\sum_{k=1}^{n} H_{k}(u,v)}{\sum_{k=1}^N H_{k}(u,v)}\le \frac{n}{N},\qquad n\le N,
\end{equation*}
which is the desired result.
\end{Proof}

Lemmas~\ref{lemma1}, \ref{lemma2} and \ref{lemma_fresh} enable us to
obtain lower pointwise estimates for the permutations $p_{K_t}$ via the
permutations $p_{\kappa_n}$ for some $t$. To do so, we will use~(\ref{phi_K1_K2}) and~(\ref{connection_kernels_n_and_m})
to estimate the coefficients of the quadratic polynomial~(\ref{parabola_t}). Let us start
with the case $n< N\le 2n$.
\begin{lemma}
\label{lemma_new1} Given $K_t$ of the form $(\ref{kernels_n_N})$
with $n< N\le 2n$, if
$$
t\in \mathbb{R}\setminus \left[-\tfrac{1}{2}\left(3+\sqrt{9-4\tfrac{N}{n}}\right);2-\tfrac{N}{n}\right],
$$
then $p_{K_t}\ge C(t)\cdot p_{\kappa_n}$ with some $C(t)>0$.
\end{lemma}
\begin{Proof}
To get the required estimate, we first look at the expression for
$\tau_2$ in (\ref{gamma2}) for our case. Since $n< N \le 2n$, from
(\ref{p_repr}) and (\ref{gamma2}) we immediately get
$$
\tau_2=\Lambda^{2(N-n)}\cdot p_{\kappa_{2n-N}},\qquad 0\le 2n-N \le n-1,
$$
with $\tau_2\equiv 0$ if $N=2n$. Consequently, by (\ref{parabola_t})
 and (\ref{phi_K1_K2}),
\begin{equation}
\label{p_{K_t}}
p_{K_t}=p_{\kappa_N}+(\tau_1 p_{\kappa_n}-\Lambda^{2(N-n)}p_{\kappa_{2n-N}})\,t+p_{\kappa_n}t^2.
\end{equation}

Now we show that the right hand side of (\ref{p_{K_t}}) for $t$ mentioned in the lemma is bounded
from below by $p_{\kappa_n}$, multiplied by a positive constant, depending only on $t$.

Applying the inequality (\ref{connection_kernels_n_and_m}) to
$p_{\kappa_N}$ and $p_{\kappa_{2n-N}}$ in (\ref{p_{K_t}}) for $t\ge
0$ gives
\begin{equation}
\label{p_K>f p_kappa} p_{K_t}\ge
\left(\tfrac{N}{n}\;\Lambda^{2(N-n)}+(\tau_1-2+\tfrac{N}{n})\;t+t^2\right)\cdot
p_{\kappa_n} =f(\xi_1,\xi_2,\xi_3)\cdot p_{\kappa_n},
\end{equation}
where
$\xi_j:=\lambda_j^{2(N-n)}\in[0,1]$, $j=1,2,3$, and
\begin{equation}
\label{f_n_N}
f(\xi_1,\xi_2,\xi_3):=\tfrac{N}{n}\;\xi_1\xi_2\xi_3+(\xi_1+\xi_2+\xi_3-2+\tfrac{N}{n})\;t+t^2.
\end{equation}
Analysis of $\partial f/\partial\xi_j$ shows that $f$ is
non-decreasing for $t\ge 0$ with respect to each $\xi_j\in[0,1]$.
Consequently,
$$
f(\xi_1,\xi_2,\xi_3)\ge f(0,0,0)=t(t-2+\tfrac{N}{n}),
$$
which is strictly positive for $t>2-\tfrac{N}{n}\ge 0$.

For $t\le 0$ we apply~(\ref{connection_kernels_n_and_m}) to $p_{\kappa_N}$ and
use that $p_{\kappa_{2n-N}}\ge 0$ (see
(\ref{positive_permut_kappa_m})). This yields
\begin{equation}
\label{p_K>F p_kappa} p_{K_t}\ge
\left(\tfrac{N}{n}\;\Lambda^{2(N-n)}+\tau_1\;t+t^2\right)\cdot
p_{\kappa_n} =F(\xi_1,\xi_2,\xi_3)\cdot p_{\kappa_n},
\end{equation}
where the function
\begin{equation}
\label{F_n_N}
F(\xi_1,\xi_2,\xi_3):=\tfrac{N}{n}\;\xi_1\xi_2\xi_3+(\xi_1+\xi_2+\xi_3)\;t+t^2
\end{equation}
is non-increasing for $t\le -\tfrac{N}{n}$ with respect to each
$\xi_j\in[0,1]$. Consequently,
$$
F(\xi_1,\xi_2,\xi_3)\ge F(1,1,1)=\tfrac{N}{n}+3\;t+t^2,
$$
where the latter expression is positive for
$t<-\tfrac{1}{2}\left(3+\sqrt{9-4\tfrac{N}{n}}\right)\le-\tfrac{N}{n}$.
\end{Proof}

Now let $N>2n$. Note that the following lemma
coincides with the previous one if we put $N=2n$.
\begin{lemma}
\label{lemma_new2} Given $K_t$ of the form $(\ref{kernels_n_N})$
with $N>2n$, if
$$
t\in \mathbb{R}\setminus \left[-\tfrac{1}{2}\left(3+\rho_{n,N}+\sqrt{(3+\rho_{n,N})^{2}-4\tfrac{N}{n}}\right); \rho_{n,N}\right],\quad \rho_{n,N}=\left(\tfrac{N}{n}-2\right)\sqrt{N-2n},
$$
then $p_{K_t}\ge C(t)\cdot p_{\kappa_n}$ with some $C(t)>0$.
\end{lemma}
\begin{Proof}
We will again estimate the coefficients of the polynomial~(\ref{parabola_t}) in terms of~$p_{\kappa_n}$. At first, we will estimate $|\phi_{n,N}|$. By (\ref{phi_K1_K2}), this will only need to estimate $|\tau_2|$.

As we have already mentioned before Lemma~\ref{lemma1}, the permutations $p_{K_t}$ and $p_{\kappa_n}$ are invariant under translations. Therefore  we can assume without loss of generality that all triangles $(0,u,v)$ that we consider  belong to
the half plane $\Re z\ge 0$. This will be necessary in the further analysis of angles of these triangles.

From now on,  we use the following
notation additionally to (\ref{notation}):
\begin{equation}
\label{notation-2}
\sin \alpha_j:=\lambda_j,\qquad \lambda_j\in [-1;1],\qquad j=1,2,3.
\end{equation}
We also suppose that $\lambda_j^2$ are pairwise distinct. One can get
the other case by passage to a limit below. For the geometrical interpretation of~$\alpha_j$ see
 Figures~\ref{fig3} and~\ref{fig4}.

Now we aim to represent $\tau_2$ from (\ref{gamma2}) in terms of the
curvature written in the form (\ref{curvature_geometry}). For this
purpose we will segregate the area squared $S(0,u,v)^2$ in the
numerator and $|u|^2|v|^2|u-v|^2$ in the denominator of $\tau_2$.
First, from (\ref{gamma2}), taking into account (\ref{notation}), we
obtain
\begin{align*}
\tau_2 &=\frac{\Lambda^{2n-1}}{|u||v||u-v|}
\frac{(\lambda_1\lambda_2)^{2(2n-N)-1}|u-v|+(\lambda_1\lambda_3)^{2(2n-N)-1}|v|
-(\lambda_2\lambda_3)^{2(2n-N)-1}|u|}{\Lambda^{2(2n-N)-1}}\\
&=\frac{\Lambda^{2n-1}}{|u||v||u-v|}\left(
\lambda_3^{2(N-2n)}\Re(u-v)
+\lambda_2^{2(N-2n)}\Re v
-\lambda_1^{2(N-2n)}\Re u\right)\\
&=\frac{\Lambda^{2n-1}}{|u||v||u-v|} \left( \Re
u\left(\lambda_3^{2(N-2n)}-\lambda_1^{2(N-2n)}\right) -\Re
v\left(\lambda_3^{2(N-2n)}-\lambda_2^{2(N-2n)}\right) \right)\\
&=\frac{\Lambda^{2n-1}}{|u||v||u-v|} \left(
\lambda_1|u|\left(\lambda_3^{2}-\lambda_1^{2}\right)A_1(u,v)
-\lambda_2|v|\left(\lambda_3^{2}-\lambda_2^{2}\right)A_2(u,v)
\right),
\end{align*}
where
\begin{equation}
\label{A1A2}
A_1(u,v):=\frac{\lambda_3^{2(N-2n)}-\lambda_1^{2(N-2n)}}{\lambda_3^{2}-\lambda_1^{2}} \quad \text{and}\quad
A_2(u,v):=\frac{\lambda_3^{2(N-2n)}-\lambda_2^{2(N-2n)}}{\lambda_3^{2}-\lambda_2^{2}}.
\end{equation}

Finally, we can rewrite $\tau_2$ as
\begin{equation}
\label{tau_2_representation}
\tau_2=\frac{\Lambda^{2n-1}}{|u|^2|v|^2|u-v|^2}\cdot A(u,v),
\end{equation}
where
$
A(u,v):=|u||v||u-v|\left(\lambda_1
\left(\lambda_3^{2}-\lambda_1^{2}\right)|u|A_1(u,v)-\lambda_2
\left(\lambda_3^{2}-\lambda_2^{2}\right)|v|A_2(u,v)\right).
$

By (\ref{notation-2}) and the formulas for the sum of sines and the sine of a double angle,
\begin{align*}
\lambda_3^{2}-\lambda_1^{2}
&=(\sin\alpha_3+\sin\alpha_1)(\sin\alpha_3-\sin\alpha_1)\\
&=2\sin\frac{\alpha_3+\alpha_1}{2}\cos\frac{\alpha_3-\alpha_1}{2}
\cdot
2\sin\frac{\alpha_3-\alpha_1}{2}\cos\frac{\alpha_3+\alpha_1}{2}\\
&= \sin(\alpha_3+\alpha_1)\sin(\alpha_3-\alpha_1).
\end{align*}
Analogously, $\lambda_3^{2}-\lambda_2^{2}=\sin(\alpha_3+\alpha_2)\sin(\alpha_3-\alpha_2)$.
Thus
\begin{equation}
\label{A}
\begin{array}{rl}
  A(u,v)= & |u|^2|v||u-v| \sin(\alpha_3+\alpha_1)\sin(\alpha_3-\alpha_1)\lambda_1 A_1(u,v)\\
           &
           -|u||v|^2|u-v|\sin(\alpha_3+\alpha_2)\sin(\alpha_3-\alpha_2)\lambda_2 A_2(u,v).
\end{array}
\end{equation}

Now let us see how one can calculate the angles $\angle(u,0,v)$,
$\angle(0,u,v)$ and $\angle(0,v,u)$ of the triangle $(0,u,v)$, using the angles $\alpha_j$,
$j=1,2,3$. Recall that the triangles $(0,u,v)$ belong to the half plane $\Re z\ge 0$. Thus only two cases are possible:
\begin{enumerate}[1.]
  \item The vertexes $u$ and $v$ both lie in the same (first or forth) quarter of the plane.
  \item The vertexes $u$ and $v$ lie in different quarters of the plane.
\end{enumerate}
\begin{figure}
  \includegraphics[width=15cm]{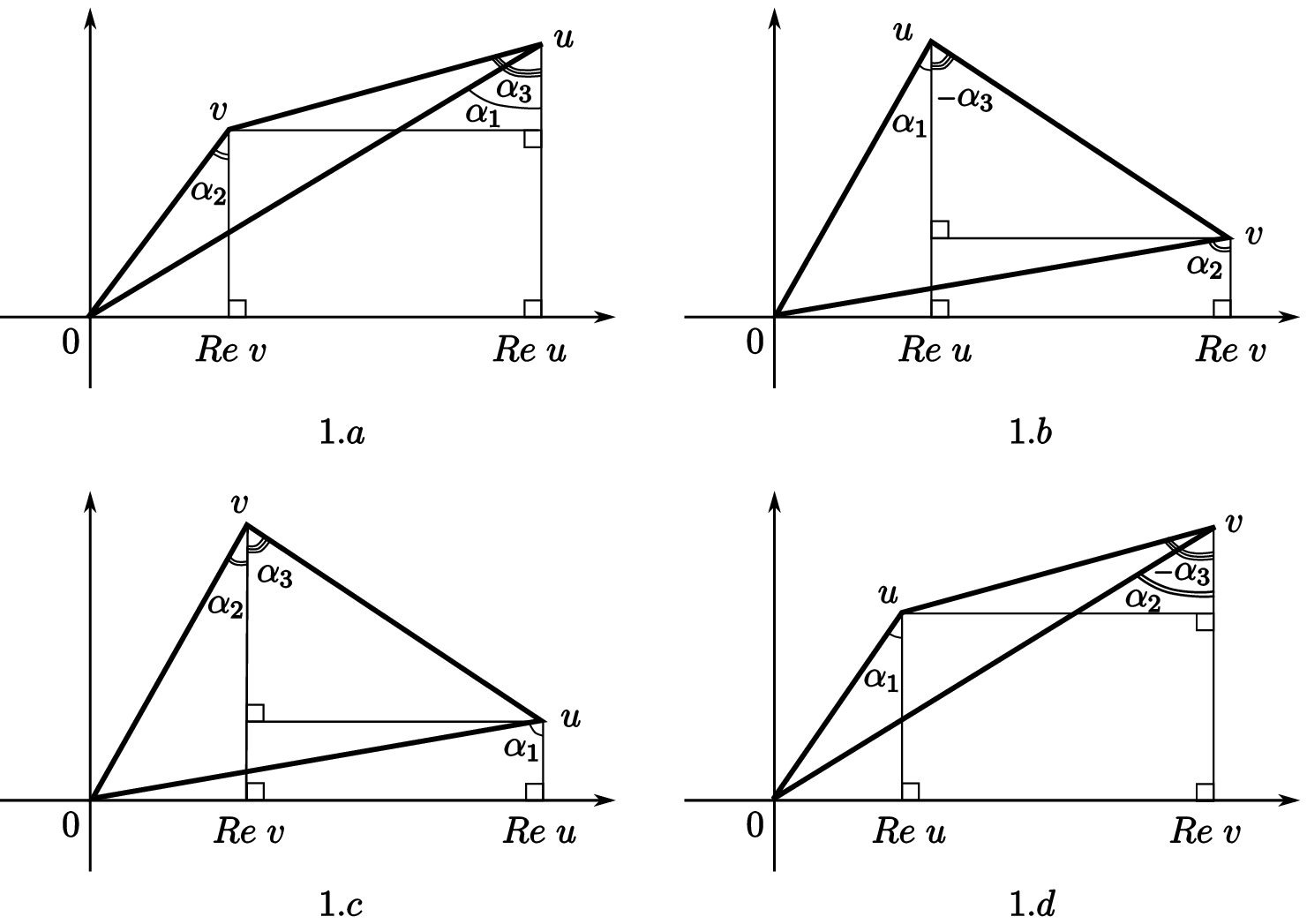}\\
 \caption{}
  \label{fig3}
\end{figure}
\begin{figure}
  \includegraphics[width=15cm]{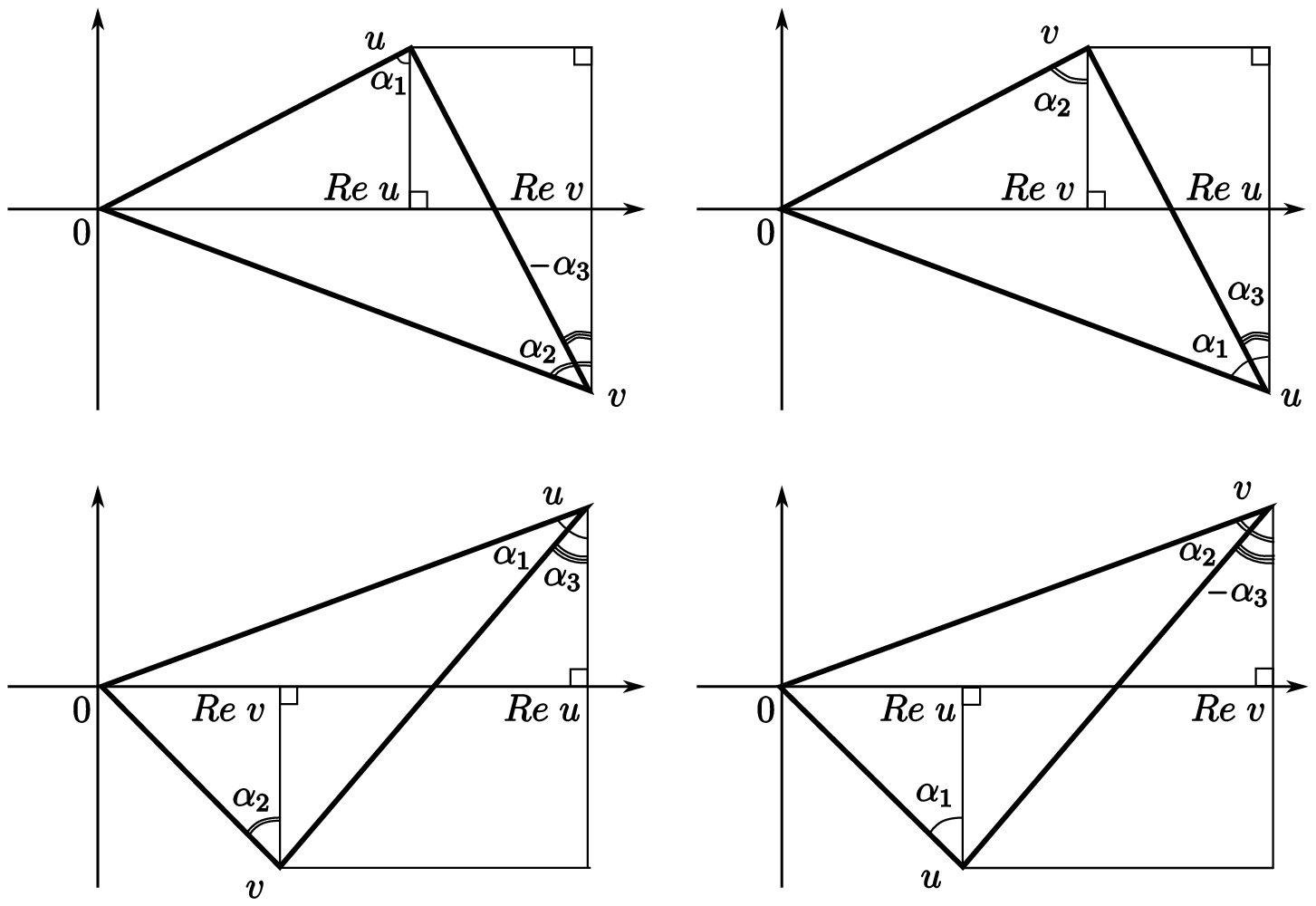}\\
 \caption{}
  \label{fig4}
\end{figure}

One can check that four options are realizable in the case~1  (see
the examples in Figure \ref{fig3}; several other situations are
possible but they produce the same cases):
\begin{enumerate}[{1}a.]
\item $\angle(u,0,v)=\phantom{-(}\alpha_1-\alpha_2\phantom{)}$, $\angle(0,u,v)=-(\alpha_1-\alpha_3)$, $\angle(0,v,u)=\pi+(\alpha_2-\alpha_3)$;
\item $\angle(u,0,v)=-(\alpha_1-\alpha_2)$, $\angle(0,u,v)=\phantom{-(}\alpha_1-\alpha_3\phantom{)}$, $\angle(0,v,u)=\pi-(\alpha_2-\alpha_3)$;
\item $\angle(u,0,v)=\phantom{-(}\alpha_1-\alpha_2\phantom{)}$, $\angle(0,u,v)=\pi-(\alpha_1+\alpha_3)$, $\angle(0,v,u)=\phantom{-(}\alpha_2+\alpha_3\phantom{)}$;
\item $\angle(u,0,v)=-(\alpha_1-\alpha_2)$, $\angle(0,u,v)=\pi+(\alpha_1+\alpha_3)$, $\angle(0,v,u)=-(\alpha_2+\alpha_3)$.
\end{enumerate}

In the case~2 (see Figure~\ref{fig4}) one always has
$$
\angle(u,0,v)=\pi-(\alpha_1+\alpha_2), \quad
\angle(0,u,v)=\alpha_1-\alpha_3,\quad
\angle(0,v,u)=\alpha_2+\alpha_3.
$$

Consequently, taking into account the formulas
$$
S(0,u,v)=\tfrac{1}{2}|u||v|\sin\angle(u,0,v)=\tfrac{1}{2}|u||u-v|\sin\angle(0,u,v)=\tfrac{1}{2}|v||u-v|\sin
\angle(0,v,u),
$$
we conclude from (\ref{A}) that,
\begin{itemize}
  \item in the cases  1.a and 1.b:
  \begin{align*}
A(u,v)&=|u||v|
    \sin(\pm(\alpha_2-\alpha_1))
  |u||u-v|
    \sin(\pm(\alpha_1-\alpha_3))
  \frac{\sin(\alpha_3+\alpha_1)}{\sin(\alpha_1-\alpha_2)}\lambda_1 A_1(u,v)\\
           &\quad -|u||v|\sin(\pm(\alpha_2-\alpha_1))
           |v||u-v|\sin(\pi\pm (\alpha_3-\alpha_2))
           \frac{\sin(\alpha_3+\alpha_2)}{\sin(\alpha_1-\alpha_2)}\lambda_2
           A_2(u,v)\\
  &= 4S(0,u,v)^2\,
  \frac{
  \sin(\alpha_3+\alpha_1)\lambda_1 A_1(u,v)-\sin(\alpha_3+\alpha_2)\lambda_2 A_2(u,v)
  }{\sin(\alpha_1-\alpha_2)};\\
\end{align*}
  \item in the cases 1.c and 1.d:
  \begin{align*}
A(u,v)&= |u||v|
    \sin(\pm(\alpha_1-\alpha_2))
  |u||u-v|
    \sin(\pi\mp(\alpha_1+\alpha_3))
  \frac{\sin(\alpha_3-\alpha_1)}{\sin(\alpha_1-\alpha_2)}\lambda_1 A_1(u,v)\\
          &\quad -|u||v|\sin(\pm(\alpha_1-\alpha_2))
           |v||u-v|\sin(\pm (\alpha_2+\alpha_3))
           \frac{\sin(\alpha_3-\alpha_2)}{\sin(\alpha_1-\alpha_2)}\lambda_2
           A_2(u,v)\\
  &= 4S(0,u,v)^2\,
  \frac{
  \sin(\alpha_3-\alpha_1)\lambda_1 A_1(u,v)-\sin(\alpha_3-\alpha_2)\lambda_2 A_2(u,v)
  }{\sin(\alpha_1-\alpha_2)};
\end{align*}
\item in the case 2:
  \begin{align*}
   A(u,v)&=-|u||v|
    \sin(\alpha_1+\alpha_2)
  |u||u-v|
    \sin(\alpha_1-\alpha_3)
  \frac{\sin(\alpha_3+\alpha_1)}{\sin(\alpha_1+\alpha_2)}\lambda_1 A_1(u,v)\\
           &\quad-|u||v|\sin(\alpha_1+\alpha_2)|v||u-v|\sin(\alpha_2+\alpha_3)
           \frac{\sin(\alpha_3-\alpha_2)}{\sin(\alpha_1+\alpha_2)}\lambda_2
           A_2(u,v)\\
  &= -4S(0,u,v)^2\,
  \frac{\sin(\alpha_3+\alpha_1)\lambda_1 A_1(u,v)+\sin(\alpha_3-\alpha_2)\lambda_2 A_2(u,v)}{\sin(\alpha_1+\alpha_2)}.
\end{align*}
\end{itemize}

Note that the substitutions $\alpha_1\mapsto -\alpha_1$,
$\alpha_2\mapsto -\alpha_2$ ($\lambda_1\mapsto -\lambda_1$,
$\lambda_2\mapsto -\lambda_2$) in the expression for $A(u,v)$ for
the case 1.a and 1.b give $A(u,v)$ in the cases 1.c and 1.d.
Moreover, the substitution  $\alpha_2\mapsto -\alpha_2$
($\lambda_2\mapsto -\lambda_2$) in $A(u,v)$ for the case 1.a and 1.b
gives $-A(u,v)$ in the case 2.
In what follows, this allows us to consider only one expression for $A(u,v)$, say, the one corresponding to 1.a and 1.b, instead of the three. This reduction will not affect the final result. By this reason, let
$$
A(u,v)=4S(0,u,v)^2\cdot \frac{V(u,v)}{\sin(\alpha_1-\alpha_2)},
$$
where
\begin{equation}
\label{V}
V(u,v):=\sin(\alpha_3+\alpha_1)\lambda_1 A_1(u,v)-
\sin(\alpha_3+\alpha_2) \lambda_2 A_2(u,v).
\end{equation}
From this and (\ref{tau_2_representation}) by the formula
(\ref{curvature_geometry}), connecting the curvature $c(0,u,v)$ and
the area $S(0,u,v)$, we get
$$
\tau_2=\frac{4S(0,u,v)^2}{|u|^2|v|^2|u-v|^2} \cdot \Lambda^{2n-1}
\cdot \frac{V(u,v)}{\sin(\alpha_1-\alpha_2)}= \tfrac{1}{4}c(0,u,v)^2
\cdot \Lambda^{2n-1} \cdot \frac{V(u,v)}{\sin(\alpha_1-\alpha_2)}.
$$
Note that $\tfrac{1}{4}c(0,u,v)^2=p_{\kappa_1}(0,u,v)$ by
(\ref{permutation_real_part}). Consequently, the inequality
(\ref{connection_kernels_n_and_m}) and the fact that $|\Lambda|\le
1$ yield
\begin{equation}
\label{|tau2|}
|\tau_2|=
 n\,\Lambda^{2(n-1)}p_{\kappa_1} \cdot\frac{|\Lambda|}{n}\cdot \frac{\left|V(u,v)\right|}{\left|\sin(\alpha_1-\alpha_2)\right|}
\le \frac{p_{\kappa_n}}{n}\cdot
\frac{|V(u,v)|}{|\sin(\alpha_1-\alpha_2)|}.
\end{equation}
Now we want to show that $|V(u,v)|\le \textrm{const} \cdot
|\sin(\alpha_1-\alpha_2)|$. If we rewrite $A_1(u,v)$ and $A_2(u,v)$,
defined in (\ref{A1A2}), using the formula
$$
\frac{a^m-b^m}{a-b}=\sum_{v=0}^{m-1}a^{m-1-v}b^v, \qquad m\in \mathbb{N}^+,
$$
for $m:=N-2n\ge 1$, then (\ref{V}) takes the form
\begin{align*}
V(u,v)= \sum_{v=0}^{m-1}\lambda_3^{2(m-1-v)}
\left(\sin(\alpha_3+\alpha_1) \cdot
\lambda_1^{2v+1}-\sin(\alpha_3+\alpha_2)\cdot
\lambda_2^{2v+1}\right).
\end{align*}
Now we substitute $\lambda_j=\sin\alpha_j$, $j=1,2$, by (\ref{notation-2}) and apply the
well-known formula
$$
(\sin \theta)^{2v+1}=\frac{1}{2^{2v}}\sum_{k=0}^v
(-1)^{v-k}\binom{2v+1}{k}\sin(2v+1-2k)\theta.
$$
This leads to the following representation:
\begin{equation}
\label{|V|}
V(u,v)=\sum_{v=0}^{m-1}\lambda_3^{2(m-1-v)}\frac{1}{2^{2v}}\sum_{k=0}^v(-1)^{v-k}\binom{2v+1}{k}B_{v,k}(\alpha_1,\alpha_2,\alpha_3),
\end{equation}
where
\begin{align*}
B_{v,k}&(\alpha_1,\alpha_2,\alpha_3)\\
&:=\sin(\alpha_3+\alpha_1)\sin(2v+1-2k)\alpha_1-\sin(\alpha_3+\alpha_2)
\sin(2v+1-2k)\alpha_2.
\end{align*}
By the formulas for the product of sines and the difference of
cosines we obtain
\begin{align*}
B_{v,k}&(\alpha_1,\alpha_2,\alpha_3)\\
=&\tfrac{1}{2}\left(\cos(\alpha_3-2(v-k)\alpha_1)-\cos(\alpha_3+2(v-k+1)\alpha_1)\right.\\
 &\left.-\cos(\alpha_3-2(v-k)\alpha_2)+\cos(\alpha_3+2(v-k+1)\alpha_2)\right)\\
=&\sin(\alpha_3-(v-k)(\alpha_1+\alpha_2))\sin((v-k)(\alpha_1-\alpha_2))\\
            &+\sin(\alpha_3+(v-k+1)(\alpha_1+\alpha_2))\sin((v-k+1)(\alpha_1-\alpha_2)).
\end{align*}
Since $|\sin r x|\le r|\sin x|$, $r\ge 0$, it follows that
$$
|B_{v,k}(\alpha_1,\alpha_2,\alpha_3)|\le
(2v-2k+1)|\sin(\alpha_1-\alpha_2)|.
$$
This and the obvious estimate of $|V(u,v)|$ from (\ref{|V|}) yield
$$
\frac{|V(u,v)|}{|\sin(\alpha_2-\alpha_1)|}\le\sigma(m):=\sum_{v=0}^{m-1}\frac{1}{2^{2v}}\sum_{k=0}^v\binom{2v+1}{k}(2v-2k+1).
$$
One can check by successive use of the
formulas (4.2.1.6), (4.2.2.13) and (4.2.3.19) from \cite[\S4.2]{Prudnikov} that
$$
\sigma(m)=\frac{4m^2-1}{3\cdot 4^{m-1}}\binom{2m-2}{m-1}.
$$
Moreover, it can by easily proved by induction that
$$
\sigma(m)\le m^{3/2},\qquad m\in \mathbb{N}.
$$
Since $m=N-2n\ge 1$, (\ref{|tau2|}) yields
\begin{equation}
\label{|tau2est|}
 |\tau_2|\le \rho_{n,N}\cdot p_{\kappa_n}, \qquad
\rho_{n,N}=\frac{(N-2n)^{3/2}}{n}=\left(\tfrac{N}{n}-2\right)\sqrt{N-2n},\qquad
N> 2n.
\end{equation}

Now we come back to the representation (\ref{parabola_t}) from Lemma
\ref{lemma1} and estimation of its terms. By (\ref{phi_K1_K2}),  (\ref{connection_kernels_n_and_m}) and
(\ref{|tau2est|}), we deduce for $t\ge 0$ that
\begin{equation}
\label{p_K>g p_kappa}
 p_{K_t}\ge
\left(\tfrac{N}{n}\;\Lambda^{2(N-n)}+(\tau_1-\rho_{n,N})\;t+t^2\right)\cdot
p_{\kappa_n} =g(\xi_1,\xi_2,\xi_3)\cdot p_{\kappa_n},
\end{equation}
where $\xi_j=\lambda_j^{2(N-n)}\in[0,1]$ as in the proof of the previous lemma, and
\begin{equation}
\label{phi_n_N}
g(\xi_1,\xi_2,\xi_3):=\tfrac{N}{n}\;\xi_1\xi_2\xi_3+(\xi_1+\xi_2+\xi_3-\rho_{n,N})\;t+t^2.
\end{equation}
The function $g$ is non-decreasing for $t\ge 0$ with respect to each
$\xi_j\in[0,1]$ , hence for ${t>\rho_{n,N}>0}$ we obtain the inequality
$$
g(\xi_1,\xi_2,\xi_3)\ge g(0,0,0)=t(t-\rho_{n,N})>0.
$$

For $t\le 0$ we have
\begin{equation}
\label{p_K>G p_kappa} p_{K_t}\ge
\left(\tfrac{N}{n}\;\Lambda^{2(N-n)}+(\tau_1+\rho_{n,N})\;t+t^2\right)\cdot
p_{\kappa_n}=G(\xi_1,\xi_2,\xi_3)\cdot p_{\kappa_n},
\end{equation}
where the function
\begin{equation}
\label{Phi_n_N}
G(\xi_1,\xi_2,\xi_3):=\tfrac{N}{n}\;\xi_1\xi_2\xi_3+(\xi_1+\xi_2+\xi_3+\rho_{n,N})\;t+t^2
\end{equation}
is non-increasing for $t\le-\tfrac{N}{n}$ with respect to each
$\xi_j\in[0,1]$ and therefore
$$
G(\xi_1,\xi_2,\xi_3)\ge G(1,1,1)=\tfrac{N}{n}+(3+\rho_{n,N})\;t+t^2.
$$
The roots of the latter quadratic polynomial are
$$
-\tfrac{1}{2}\left(3+\rho_{n,N}\pm\sqrt{(3+\rho_{n,N})^{2}-4\tfrac{N}{n}}\right),
$$
so it has only positive values if $t<-\tfrac{1}{2}\left(3+\rho_{n,N}+\sqrt{(3+\rho_{n,N})^{2}-4\tfrac{N}{n}}\right)\le -\tfrac{N}{n}$.
\end{Proof}

\medskip

Note that Lemmas~\ref{lemma_new1} and~\ref{lemma_new2} give
Theorem~\ref{new_theorem1} by continuity. For the proof of
Theorem~\ref{new_theorem2}, additionally to
Theorem~\ref{new_theorem1}, we will also need lower estimates of
$p_{K_t}$ for~$t$, which are the end points of the intervals excluded in (\ref{estimates_t_N<2n}) and (\ref{estimates_t_N>2n}) from the real line.  And to obtain these estimates, we first introduce additional notation.

Given two distinct points $z,w\in \mathbb{C}$, we denote by
$L_{z,w}$ the line passing through $z$ and $w$. Given three pairwise
distinct points $z_1,z_2,z_3\in \mathbb{C}$, we denote by
$\angle(z_1,z_2,z_3)$ the smallest angle formed by the lines
$L_{z_1,z_2}$ and $L_{z_1,z_3}$. This angle belongs to $[0;\pi/2]$.
If $L$ and $L'$ are lines, then $\angle(L,L')$ is the smallest angle
between them. This angle belongs to $[0;\pi/2]$, too. Also,
$\theta_V(L):=\angle(L,V)$ and $\theta_H(L):=\angle(L,H)$, where $V$
and $H$ are the vertical and horizontal lines, correspondingly.
Furthermore, for a fixed constant $\tau\ge 1$, we set
\begin{equation}
\label{O_tau}
\mathcal{O}_\tau=\left\{(z_1,z_2,z_3):\frac{|z_i-z_j|}{|z_i-z_k|}\le\tau
\text{ for pairwise distinct } i,j,k\in \{1,2,3\}\right\},
\end{equation}
so that all the triangles with vertexes $z_1$, $z_2$ and $z_3$ in $\mathcal{O}_\tau$ have comparable sides.

Given $\alpha_0\in(0,\pi/2)$ and $(z_1,z_2,z_3)$, in what follows we will use the conditions
\begin{equation}
\label{delta1}
\theta_{V}(L_{z_1,z_2})+\theta_{V}(L_{z_2,z_3})+\theta_{V}(L_{z_1,z_3})\ge
\alpha_0
\end{equation}
and
\begin{equation}
\label{delta2}
\theta_{H}(L_{z_1,z_2})+\theta_{H}(L_{z_2,z_3})+\theta_{H}(L_{z_1,z_3})\ge
\alpha_0.
\end{equation}
Note that (\ref{delta1}) and (\ref{delta2}) can be  correspondingly replaced by the conditions
$$
\theta_{H}(L_{z_1,z_2})+\theta_{H}(L_{z_2,z_3})+\theta_{H}(L_{z_1,z_3})\le
\tfrac{3}{2}\pi-\alpha_0
$$
and
$$
\theta_{V}(L_{z_1,z_2})+\theta_{V}(L_{z_2,z_3})+\theta_{V}(L_{z_1,z_3})\le
\tfrac{3}{2}\pi-\alpha_0.
$$

To obtain the desired result, we first prove several geometrical lemmas.
\begin{lemma}
\label{lemma_ends1} Fix $\alpha_0\in(0,\pi/2)$. Given $(0,u,v)\in
\mathcal{O}_\tau$, if the condition $(\ref{delta1})$ is satisfied,
then
$$
\tau_1(0,u,v)=\lambda_1^{2(N-n)}+\lambda_2^{2(N-n)}+\lambda_3^{2(N-n)}\ge C_1(\alpha_0)>0.
$$
\end{lemma}
\begin{Proof}
Clearly,
$$
\lambda_1^2=\sin^2\theta_{V}(L_{0,u}), \quad \lambda_2^2=\sin^2\theta_{V}(L_{0,v}), \quad \lambda_3^2=\sin^2\theta_{V}(L_{u,v}).
$$
Moreover, from (\ref{delta1}) it follows that at least one of the angles $\theta_{V}(L_{0,u})$, $\theta_{V}(L_{0,v})$, $\theta_{V}(L_{u,v})$ is not less than $\alpha_0/3$. Thus
$\tau_1\ge (\sin\tfrac{\alpha_0}{3})^{2(N-n)}$.
\end{Proof}
\begin{lemma}
\label{lemma_ends2} Fix $\alpha_0\in(0,\pi/2)$. Given $(0,u,v)\in
\mathcal{O}_\tau$, if the condition $(\ref{delta2})$ is satisfied,
then
$$
\Upsilon(0,u,v):=2+(\lambda_1\lambda_2\lambda_3)^{2(N-n)}-(\lambda_1^{2(N-n)}+\lambda_2^{2(N-n)}+\lambda_3^{2(N-n)})\ge
C_2(\alpha_0,\tau)>0.
$$
\end{lemma}
\begin{Proof}
First we note that
$$
\Upsilon(0,u,v)\ge 2+\lambda_1^2\lambda_2^2
\lambda_3^2-(\lambda_1^2+\lambda_2^2+\lambda_3^2)
$$
as the function $2+\xi_1\xi_2\xi_3-(\xi_1+\xi_2+\xi_3)$ is
non-increasing with respect to each $\xi_j\in [0,1]$, $j=1,2,3$, and
$\lambda_j^2\ge \lambda_j^{2(N-n)}$ for $N> n\ge 1$.

In addition, we have
$$
 \lambda_1^2=1-\sin^2\theta_{H}(L_{0,u}), \quad \lambda_2^2=1-\sin^2\theta_{H}(L_{0,v}), \quad \lambda_3^2=1-\sin^2\theta_{H}(L_{u,v}),
$$
and hence
\begin{align*}
\Upsilon(0,u,v)&\ge \sin^2\theta_{H}(L_{0,u})\sin^2\theta_{H}(L_{0,v})+\sin^2\theta_{H}(L_{0,u})\sin^2\theta_{H}(L_{u,v})\\
&\qquad+\sin^2\theta_{H}(L_{0,v})\sin^2\theta_{H}(L_{u,v})-\sin^2\theta_{H}(L_{0,u})\sin^2\theta_{H}(L_{0,v})\sin^2\theta_{H}(L_{u,v})\\
&\ge \tfrac{2}{3}
\left(\sin^2\theta_{H}(L_{0,v})\sin^2\theta_{H}(L_{0,u})+\sin^2\theta_{H}(L_{0,u})\sin^2\theta_{H}(L_{u,v})\right.\\
&\qquad\left.+\sin^2\theta_{H}(L_{0,v})\sin^2\theta_{H}(L_{u,v})\right).
\end{align*}
Consider a triangle $(0,u,v)\in \mathcal{O}_\tau$ such that (\ref{delta2}) is satisfied. Fix
some $\varepsilon\in(0;\alpha_0/3)$. Two cases are possible:
\begin{enumerate}
  \item amongst $\theta_{H}(L_{0,u})$, $\theta_{H}(L_{0,v})$,
$\theta_{H}(L_{u,v})$, there exists a pair of angles, each being
greater than $\varepsilon$ and then it is easily seen that
$\Upsilon(0,u,v)\ge \tfrac{2}{3}\sin^4\varepsilon$;
  \item amongst those, there exists no pair of angles, each being
greater than $\varepsilon$.
\end{enumerate}

Let us consider the second case in detail (see Figure~\ref{fig2}).
It is clear that at least two angles amongst
$\theta_{H}(L_{0,u})$, $\theta_{H}(L_{0,v})$, $\theta_{H}(L_{u,v})$
are less than $\varepsilon$ then. In other words, two sides of the
triangle cut the horizontal line at angles less than $\varepsilon$.
We call these sides $A$ and~$B$.
\begin{figure}
  \includegraphics[width=15cm]{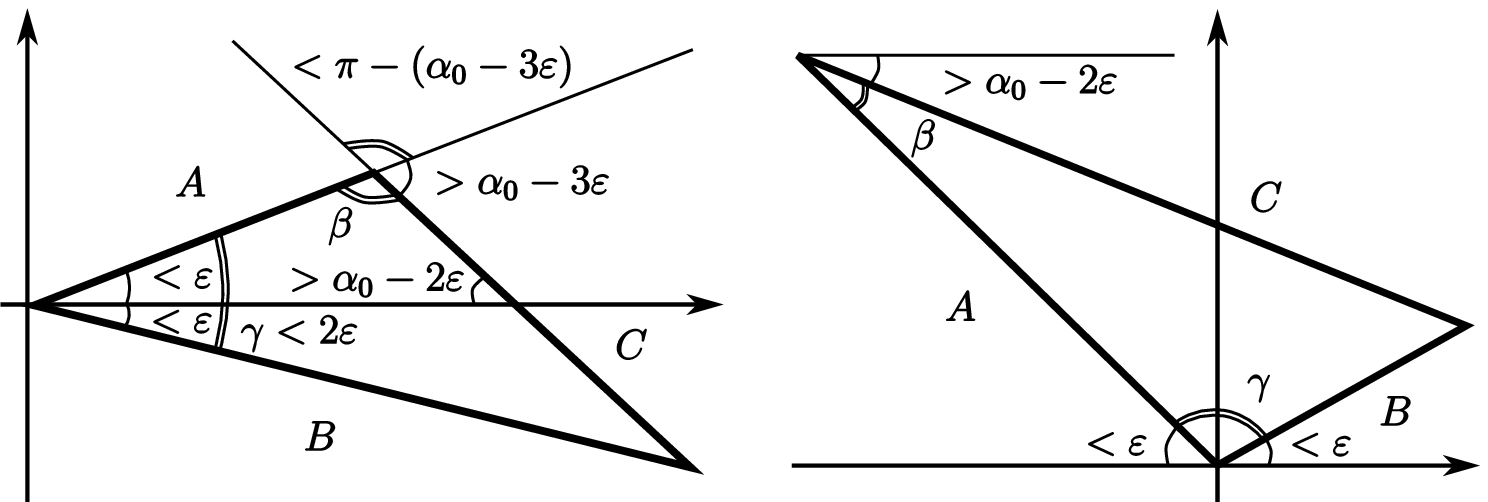}\\
\caption{}
  \label{fig2}
\end{figure}

Furthermore, let the angle $\gamma$ between $A$ and $B$ be acute;
then obviously it is smaller than $2\varepsilon$. Then the acute
angle between the third side $C$ and the horizontal line is greater
than $\alpha_0 -2\varepsilon$ and the acute angle between $A$ and $C$
is greater than $\alpha_0 -3\varepsilon$. Consequently, the obtuse
angle between $A$ and $C$ is smaller than $\pi-(\alpha_0
-3\varepsilon)$. Thus we have for the angle $\beta$ of the
triangle:
$$
\alpha_0 -3\varepsilon< \beta< \pi-(\alpha_0 -3\varepsilon).
$$
Therefore by the law of sines, the inequalities $(2/\pi)\,x\le \sin
x \le x$ for $x\in[0,\pi/2]$, and (\ref{O_tau}), we get
$$
\frac{1}{\tau^{2}}\le \frac{\textrm{length}(C)}{\textrm{length}(B)}=
\frac{\sin \gamma}{\sin \beta}<\frac{\sin
2\varepsilon}{\sin(\alpha_0 -3\varepsilon)} \le
\frac{\pi\varepsilon}{\alpha_0 -3\varepsilon} \; \Rightarrow \;
\varepsilon>\varepsilon_0(\alpha_0,\tau):=\frac{\alpha_0}{3+\pi\tau^2}.
$$

Now let the angle $\gamma$ between $A$ and $B$ be not acute (it is
greater than $\pi- 2\varepsilon$). Then for one of acute angles of
the triangle, say $\beta$, we have
$$
\beta< \varepsilon  - (\alpha_0 -2\varepsilon) =3\varepsilon-\alpha_0 <0, \qquad \varepsilon\in(0;\alpha_0/3),
$$
which is impossible.

It follows from the aforesaid that  there is a contradiction for $\varepsilon=\varepsilon_0(\alpha_0,\tau)$ in the second case and  thus   $\Upsilon(0,u,v)\ge \tfrac{2}{3}\sin^4\varepsilon_0(\alpha_0,\tau)$.
\end{Proof}

\medskip

We will also need the following result.
\begin{lemma}[Lemma 2.3 in \cite{CMPT}]
\label{lemma_CMPT_t=0} Fix $\alpha_0\in(0,\pi/2)$. Given $\kappa_n$
of the form $(\ref{kappa_n})$ and $(z_1,z_2,z_3)\in
\mathcal{O}_\tau$, if the condition $(\ref{delta1})$ is satisfied,
then
\begin{equation*}
\label{CMPT_lemma}
 p_{\kappa_n}(z_1,z_2,z_3)\ge
C_3(\alpha_0,\tau)\cdot c(z_1,z_2,z_3)^2
\end{equation*}
for some $C_3(\alpha_0,\tau)>0$.
\end{lemma}

Now we are able to obtain necessary lower pointwise estimates for $p_{K_t}$
if $t$ are the end points of the intervals  excluded in
(\ref{estimates_t_N<2n}) and (\ref{estimates_t_N>2n}) from the real line. Recall that
$\rho_{n,N}=(\tfrac{N}{n}-2)\sqrt{N-2n}$.

\begin{lemma}
\label{lemma_ends}
Fix $\alpha_0\in(0,\pi/2)$. Given $K_t$ of the form $(\ref{kernels_n_N})$ and $(z_1,z_2,z_3)\in \mathcal{O}_\tau$,
\begin{enumerate}[$($i$)$]
  \item if the condition $(\ref{delta1})$ is satisfied and $t=2-\tfrac{N}{n}$ for $n<N\le 2n$ or $t=\rho_{n,N}$ for $N\ge 2n$, or
  \item if both the conditions $(\ref{delta1})$ and $(\ref{delta2})$ are satisfied and $t=-\tfrac{1}{2}\left(3+\sqrt{9-4\tfrac{N}{n}}\right)$ for $n<N\le 2n$ or $t=-\frac{1}{2}\left(3+\rho_{n,N}+\sqrt{(3+\rho_{n,N})^2-4\frac{N}{n}}\right)$ for $N\ge 2n$,
\end{enumerate}
then
$$
p_{K_t}(z_1,z_2,z_3)\ge C(\alpha_0,\tau)\cdot c(z_1,z_2,z_3)^2
$$
for some $C(\alpha_0,\tau)>0$.
\end{lemma}
\begin{Proof}
It is enough to prove it for triples $(0,u,v)$. What is more,
the statement for $t=0$ in (\textit{i}), i.e. when $N=2n$, is just
Lemma~\ref{lemma_CMPT_t=0} and therefore we may exclude it. We also
recall the notation $\xi_j:=\lambda_j^{2(N-n)}\in[0,1]$, $j=1,2,3$.

Now let $t=2-\tfrac{N}{n}$ and $n<N< 2n$. Then for the function
given in (\ref{f_n_N}) we have
$$
f(\xi_1,\xi_2,\xi_3)=\tfrac{N}{n}\;\xi_1\xi_2\xi_3+(2-\tfrac{N}{n})(\xi_1+\xi_2+\xi_3)
\ge (2-\tfrac{N}{n})\tau_1,
$$
where $\tau_1$ is as in (\ref{gamma1}). If $t=\rho_{n,N}$ and $N>
2n$, then from (\ref{phi_n_N}) it follows that
$$
g(\xi_1,\xi_2,\xi_3)=\tfrac{N}{n}\;\xi_1\xi_2\xi_3+\rho_{n,N}(\xi_1+\xi_2+\xi_3)
\ge \tfrac{N}{n}\tau_1.
$$

Thus from the inequalities (\ref{p_K>f p_kappa}) for $n<N<2n$ and
(\ref{p_K>g p_kappa}) for $N> 2n$, both being valid for $t\ge 0$,
and Lemmas~\ref{lemma_ends1} and~\ref{lemma_CMPT_t=0} (with the
assumption $(\ref{delta1})$) we get
$$
p_{K_t}(0,u,v)\ge \tfrac{N}{n}  \tau_1(0,u,v) \cdot
p_{\kappa_n}(0,u,v)\ge \tfrac{N}{n}\,
C_1(\alpha_0)\,C_3(\alpha_0,\tau)\cdot c(0,u,v)^2,
$$
which is the required result in the case (\textit{i}).

Now consider (\textit{ii}). Set
$t=-\tfrac{1}{2}\left(3+\sqrt{9-4\tfrac{N}{n}}\right)$ and $n<N\le
2n$. Then, by (\ref{F_n_N}),
$$
F(\xi_1,\xi_2,\xi_3)=\tfrac{N}{n}(\xi_1\xi_2\xi_3-1)+\tfrac{1}{2}(3-(\xi_1+\xi_2+\xi_3))(3+\sqrt{9-4\tfrac{N}{n}}).
$$
Since $-\tfrac{1}{2}\left(3+\sqrt{9-4\tfrac{N}{n}}\right)\le -\tfrac{N}{n}$,
$$
F(\xi_1,\xi_2,\xi_3)\ge
\tfrac{N}{n}(2+\xi_1\xi_2\xi_3-(\xi_1+\xi_2+\xi_3))\ge
\tfrac{N}{n}\Upsilon.
$$

The function (\ref{Phi_n_N}) for
$$
t=t_0:=-\tfrac{1}{2}\left(3+\rho_{n,N}+\sqrt{(3+\rho_{n,N})^2-4\tfrac{N}{n}}\right)
$$
and ${N\ge 2n}$ after some simplifications takes the form
$$
G(\xi_1,\xi_2,\xi_3)=\tfrac{N}{n}(\xi_1\xi_2\xi_3-1)+(\xi_1+\xi_2+\xi_3-3)t_0,
$$
and hence $G(\xi_1,\xi_2,\xi_3)\ge \tfrac{N}{n}\Upsilon$ since
$t_0\le -\tfrac{N}{n}$ for $N\ge 2n$.

Thus  for the last two values of $t$, by the inequalities
(\ref{p_K>F p_kappa}) for $n<N\le 2n$ and (\ref{p_K>G p_kappa}) for
$N\ge 2n$, both being valid for $t\le -\tfrac{N}{n}$, and Lemma
\ref{lemma_CMPT_t=0} (with the assumption $(\ref{delta1})$), we get
$$
p_{K_t}(0,u,v)\ge \tfrac{N}{n}\,\Upsilon(0,u,v)\cdot
p_{\kappa_n}(0,u,v)\ge \tfrac{N}{n}\,\Upsilon(0,u,v)\,
C_3(\alpha_0,\tau)\cdot c(0,u,v)^2.
$$
If  $(\ref{delta2})$ is also satisfied, then by
Lemma~\ref{lemma_ends2} we obtain the desired inequality
$$
p_{K_t}(0,u,v)\ge \tfrac{N}{n}\,C_2(\alpha_0,\tau)\,
C_3(\alpha_0,\tau)\cdot c(0,u,v)^2,
$$
and we are done.
\end{Proof}

\section{Examples}
\label{section_examples}
 In this section we  present triples $(0,u,v)$ such that the
permutations $p_{K_t}(0,u,v)$ change sign for $t$ mentioned in~(\ref{t_bounds}),
namely,
 $t\in (-N/n;0)$ for all $n$ and $N$  and $t\in\left(0;\frac{2}{e^{3/2}}\frac{N}{n}\right)$ for $N\gg n$.
We use the notations of Lemma~\ref{lemma1} below. Note that by this
lemma, $p_{K_t}$ can be calculated via $p_{\kappa_m}$ for $m$, equal
to $n$ and $N$, and $\varphi_{n,N}$. To obtain $p_{\kappa_m}$ we
substitute (\ref{kappa_n}) into (\ref{p_K}); $\varphi_{n,N}$ is
calculated by (\ref{varphi}).

We first show that $p_{K_t}(0,u,v)$ is positive for any $t$ if $u=a+i$ and
$v=a-i$, where $a\in \mathbb{R}\setminus \{0\}$ is suitably chosen. By (\ref{p_K}) and taking into account that
$K_t(u)=K_t(v)$ and $K_t(u-v)\equiv 0$,
$$
p_{K_t}
=\left(\frac{a^{2N-1}}{(1+a^2)^{N}}+t\cdot \frac{a^{2n-1}}{(1+a^2)^{n}}\right)^2
=\frac{a^{4n-2}}{(1+a^2)^{2n}}\left(\frac{a^{2(N-n)}}{(1+a^2)^{N-n}}+t\right)^2,
$$
which is positive for any real $t$ if $a$ is chosen so that the expression in the latter brackets does not vanish.

Now the aim is to show that for any fixed $t$ from (\ref{t_bounds})
there exist triples $(0,u,v)$ such that $p_{K_t}(0,u,v)$ is
negative. To do so, we find families of $(0,u,v)$ such that the
quadratic polynomial
$$
\frac{p_{K_t}}{p_{\kappa_n}}=\frac{p_{\kappa_N}}{p_{\kappa_n}}+\frac{\varphi_{n,N}}{p_{\kappa_n}}t+t^2, \qquad p_{\kappa_n}>0,
$$
(with respect to $t$) has two different roots $t_1$ and $t_2$,
depending on $u$ and $v$, and thus $p_{K_t}$ (of the form
(\ref{parabola_t})) is negative for $t\in (t_1;t_2)$. In addition, we
prove that the union of the intervals $(t_1;t_2)$ when $(0,u,v)$ runs
the whole above-mentioned family is either the interval $(-N/n;0)$
or $\left(0;\frac{2}{e^{3/2}}\frac{N}{n}\right)$, indicated in
(\ref{t_bounds}).

\medskip

Let us consider the case $t<0$.
\begin{example}
\label{example1}
Set $u=-a+i$, $v=a+i$, where $a\in \mathbb{R}\setminus \{0\}$. Then
\begin{align*}
p_{\kappa_m}
=-\frac{a^{2(2m-1)}}{(a^2+1)^{2m}}+\frac{a^{2m-2}}{2(a^2+1)^{m}}
+\frac{a^{2m-2}}{2(a^2+1)^{m}}
=
\frac{a^{2m-2}\left((a^2+1)^m-a^{2m}\right)}{(a^2+1)^{2m}},
\end{align*}
where $m$ equals $n$ or $N$, and
\begin{align*}
\varphi_{n,N}=
&-\frac{a^{2N-1}}{(a^2+1)^{N}}\left(\frac{a^{2n-1}}{(a^2+1)^{n}}-\frac{1}{2a}\right)\\
&+\frac{a^{2N-1}}{(a^2+1)^{N}}\left(-\frac{a^{2n-1}}{(a^2+1)^{n}}+\frac{1}{2a}\right)
                        -\frac{1}{2a}\left(-\frac{a^{2n-1}}{(a^2+1)^{n}}-\frac{a^{2n-1}}{(a^2+1)^{n}}\right)\\
=&
\frac{a^{2n-2}\left((a^2+1)^N+(a^2+1)^na^{2(N-n)}-2a^{2N}\right)}{(a^2+1)^{N+n}}.
\end{align*}
From this by (\ref{parabola_t}) we deduce that
\begin{equation}
\label{exampl1_parab}
\frac{p_{K_t}}{p_{\kappa_n}}=
\frac{d_1(a)d_2(a)}{d_3(a)^2}
+\frac{d_1(a)+d_2(a)}{d_3(a)}t+t^2,\qquad a\neq 0,
\end{equation}
where
$$
d_1(a):=(a^2+1)^n\left(1-\frac{a^{2N}}{(a^2+1)^{N}}\right),\;\;
d_2(a):=(a^2+1)^n\left(\frac{a^{2(N-n)}}{(a^2+1)^{N-n}}-\frac{a^{2N}}{(a^2+1)^{N}}\right)
$$
and $d_3(a):=(a^2+1)^n-a^{2n}$. The polynomial (\ref{exampl1_parab})
 has two different negative roots $t_1(a)=-d_1(a)/d_3(a)$ and $t_2(a)=-d_2(a)/d_3(a)$,
 where ${d_1(a)>d_2(a)>0}$ and $d_3(a)>0$. It is easy to check that
  the roots $t_1(a)$ and $t_2(a)$ run the intervals $(-N/n;-1)$ and $(-1;0)$, correspondingly, when $a$ runs
  $(0;\infty)$. Furthermore, we see by continuity that
  $$
  \bigcup_{a\in(0;\,\infty)}(t_1(a);t_2(a))=(-N/n;0).
  $$
  This means that $p_{K_t}$ is negative for any $t$ in $(-N/n;0)$ for $a$ suitably chosen.
\end{example}

As we have already mentioned above, this example shows that
(\ref{estimates_t_N>2n}) is sharp for $N=2n$ in the sense of
Remark~\ref{remark1}. In addition, the left hand side of
(\ref{estimates_t_N>2n}) is also sharp for $N=2n+1$.

\medskip

Now we give the example for $t>0$.
\begin{example}
\label{example4} Let $n$ be fixed. Consider the triples $(0,u,v)$ such
that
\begin{equation}
\label{points}
 u=-r(1+\delta_{N}(q)i), \qquad
v=r(1-\tfrac{1}{r}+\delta_{N}(q)i),\qquad
\delta_{N}(q):=\sqrt{\frac{\ln q}{N-n}},
\end{equation}
where $r>0$ and $q\ge e$. We can also calculate
$p_{\kappa_m}$ and $\varphi_{n,N}$ for these $(0,u,v)$ using (\ref{kappa_n}),
(\ref{p_K}) and (\ref{varphi}). However, the expression of
$p_{K_t}/p_{\kappa_n}$ obtained is too big and therefore we do not
place it here. Instead, we give the following identity (the
permutations are calculated for $(0,u,v)$ as in (\ref{points})):
$$
P(t):=\lim_{r\to
\infty}\frac{p_{K_t}}{p_{\kappa_n}}=c_N(q)+b_N(q)t+t^2,
$$
where
$$
c_N(q)=\frac{N}{n\,(1+\delta^2_{N}(q))^{2(N-n)}}\cdot\frac{(2N-1)\delta^2_{N}(q)+1}{(2n-1)\delta^2_{N}(q)+1}
$$
and
$$
b_N(q)=-\frac{\left(2(N-n)^2+N-4nN+n\right)\delta^2_{N}(q)-(n+N)}{n\,(1+\delta^2_{N}(q))^{N-n}((2n-1)\delta^2_{N}(q)+1)}.
$$
Note that
$$
c_N(q)\sim \frac{(2\ln q+1)N}{q^2n},\qquad b_N(q)\sim -\frac{(2\ln
q-1)N}{qn},\qquad N\to \infty.
$$
The quadratic polynomial $P$ (with respect to $t$) has two different
 positive roots $t_1(N,q)$ and $t_2(N,q)$ if
$N$ is large enough (as the discriminant is positive). Additionally, one can check that
$$
t_1(N,q)\sim \tilde{t}_1(q):= \frac{2 \ln q+1}{q\,(2\ln q-1)},\quad
t_2(N,q) \sim \tilde{t}_{2}(N,q):=\frac{(2\ln q-1)N}{qn},\quad
N\to\infty.
$$
Taking into account the properties
$$
\tilde{t}_1(q)\to 0 \text{ as } q\to \infty \qquad \text{and} \qquad
\max_{q\in[e;\,\infty)} \frac{2\ln q-1}{q}=\frac{2}{e^{3/2}},
$$
we deduce by continuity that
$$
  \bigcup_{q\in[e;\,\infty)} (\tilde{t}_1(q);\tilde{t}_2(q))=\left(0;\frac{2}{e^{3/2}}\frac{N}{n}\right).
$$

Thus, $p_{K_t}(0,u,v)$ with $u$ and $v$ as
in~(\ref{points}) are negative for any $t$ in
$\left(0;\frac{2}{e^{3/2}}\frac{N}{n}\right)$, if $N$ (with respect
to $n$) and $r$ are large enough and $q$ is suitably chosen.
\end{example}

\section{Proof of Theorems~1 and 2}
\label{section_proof}

Recall that Lemmas~\ref{lemma_new1} and~\ref{lemma_new2} state that
if $K_t$ is of the form $(\ref{kernels_n_N})$ and
\begin{eqnarray*}
&t\in \mathbb{R}\setminus \left[-\tfrac{1}{2}\left(3+\sqrt{9-4\tfrac{N}{n}}\right); 2-\tfrac{N}{n}\right],\qquad & n<N\le 2n,\\
&t\in \mathbb{R}\setminus
\left[-\frac{1}{2}\left(3+\rho_{n,N}+\sqrt{(3+\rho_{n,N})^2-4\frac{N}{n}}\right);
\rho_{n,N}\right],\qquad &N\ge 2n,
\end{eqnarray*}
where $\rho_{n,N}=\left(\tfrac{N}{n}-2\right)\sqrt{N-2n}$, then
$$
p_{K_t}\ge C(t)\cdot p_{\kappa_n},\qquad C(t)>0.
$$

These lemmas immediately give Theorem~\ref{new_theorem1} by
continuity if we take into account the fact that
$p_{\kappa_n}(z_1,z_2,z_2)\ge 0$ for all $(z_1,z_2,z_3)\in
\mathbb{C}^3$ (see (\ref{curvature_geometry}) and
(\ref{positive_permut_kappa_m})).

\medskip

What is said from now on is related to  Theorem~\ref{new_theorem2}.

First of all we note that the latter statement of
Theorem~\ref{new_theorem2}, i.e. the one asserting that if the
operator $T_{K_t}$ is $L^2(\mathcal{H}^1\lfloor E)$-bounded, then
$E$ is rectifiable, is a corollary of the fact that the
$L^2(\mathcal{H}^1\lfloor E)$-boundedness of $T_{K_t}$ implies that
 $p_{K_t}(\mathcal{H}^1\lfloor E)<\infty$. This can be proved
by generalization of the Melnikov-Verdera identity
(\ref{Melnikov-Verdera}) for the kernels $K_t$ and permutations
$p_{K_t}$; see Lemma~3.3 in \cite{CMPT} and also Lemma 2.1 in
\cite{MMV}.

Now we come to the proof of the former statement in
Theorem~\ref{new_theorem2}.

The proof for $t$, mentioned in Lemmas~\ref{lemma_new1} and
\ref{lemma_new2} (see also the beginning of the current section), is
direct via Theorem~\ref{CMPT_theorem}, which states that if $p_{\kappa_n}(\mathcal{H}^1\lfloor E)<\infty$, then  $E$ is rectifiable. Indeed, if $p_{K_t}(\mathcal{H}^1\lfloor
E)<\infty$ for such $t$, then $p_{\kappa_n}(\mathcal{H}^1\lfloor
E)<\infty$ by the inequality $p_{K_t}\ge C(t)\cdot p_{\kappa_n}$,
$C(t)>0$, and thus the set $E$ is rectifiable.

What is left is to prove the former statement in
Theorem~\ref{new_theorem2} for
\begin{eqnarray*}
&t=-\tfrac{1}{2}\left(3+\sqrt{9-4\tfrac{N}{n}}\right),\qquad t=2-\tfrac{N}{n},\qquad & n<N\le 2n,\\
&t=-\frac{1}{2}\left(3+\rho_{n,N}+\sqrt{(3+\rho_{n,N})^2-4\frac{N}{n}}\right),\qquad
t=\rho_{n,N},\qquad&N\ge 2n.
\end{eqnarray*}
It requires some additional work and therefore for the reader's
convenience we first make several observations, which could help to
clarify the forthcoming proof.

We start with a very brief exposition of the proof of
Theorem~\ref{theorem_Leger} given in~\cite{L} (note
that one can find a modified version of the proof
from~\cite{L} in \cite[Chapter 7]{Tolsa_book} and follow it
instead). Recall that Theorem~\ref{theorem_Leger} states that, for
a Borel set $E\subset\mathbb{C}$ such that
${0<\mathcal{H}^1(E)<\infty}$, if ${c^2(\mathcal{H}^1\lfloor
E)<\infty}$, then $E$ is rectifiable. We emphasize again that it is essential in the
proof that the curvature is non-negative.

The first step is to show that there exists a compact subset $F$ of
the given set $E$ such that, among other things,
$c^2(\mathcal{H}^1\lfloor F)$ is well-controlled and can be made
very small (this is done in \cite{L} by a quite standard
uniformization procedure). Then the second and most important step
follows --- to prove that if $\mu$ is a positive Radon measure on
$\mathbb{C}$ satisfying a few special conditions (see
Proposition~\ref{proposition2} below), then there exists a Lipschitz
graph $\Gamma$ such that $\mu(\Gamma)\ge C\cdot \mu(\mathbb{C})$,
where $C$ is an absolute constant.

The problem is to choose an adequate coordinate system of
$\mathbb{C}$ and construct a Lipschitz function $A$ whose graph will
be the one needed. For this purpose, the author of \cite{L} first
defines some functions used to measure how well the $\textrm{spt}\,
\mu$ is approximated by straight lines at a given location and a
given scale. It is shown that these functions are related to the
$c^2(\mu)$ in the case when the measure $\mu$ does not degenerate
too much. These preliminary results are then used to construct the
function $A$ by stopping time arguments, which demand fine
adjustments to many parameters and thresholds. Starting with
choosing a point $x_0\in \textrm{spt}\, \mu$ and fixing an
approximating line $D_0$ (which will be the domain of the
function~$A$) such that the mean distance from $\textrm{spt}\, \mu$
to the line $D_0$ is suitably small, the author of \cite{L} comes to
cutting $\textrm{spt}\, \mu$ in four disjoint pieces $Z$, $F_1$,
$F_2$ and $F_3$ such that
$$
\textrm{spt}\, \mu=Z\cup F_1\cup F_2\cup F_3.
$$
It is shown that $Z$ is very nice for constructing the graph but the
three others admit ``bad events''. Then the goal is to prove that
these bad pieces carry only a small part of the measure $\mu$,
namely, $\mu(F_j)\le 10^{-6}\mu(\mathbb{C})$ for each $j$ and thus
$\mu(Z)\ge \frac{99}{100}\mu(\mathbb{C})$. This allows to construct
the required Lipschitz function $A: D_0 \to D_0^\perp$  such that
the set $Z$ is contained in the graph of $A$.

Coming back to the initial settings, if $\mu=\mathcal{H}^1\lfloor
F$, where $F$ is the above-mentioned subset of $E$, then there
exists a Lipschitz graph $\Gamma$ such that
$\mathcal{H}^1(\Gamma\cap F)\ge C\cdot \mathcal{H}^1(F)$. This fact
is used in the last step of the proof from \cite{L}, which is as
follows. Since $\mathcal{H}^1(E)<\infty$ by the assumptions,  the
set $E$ can be decomposed into a rectifiable and purely
unrectifiable part, i.e. $E=E_{\mathrm{rect}}+E_{\mathrm{unrect}}$.
Suppose that
\begin{equation}
\label{contradiction}
{\mathcal{H}^1(E_{\mathrm{unrect}})>0}.
\end{equation}
Then there exists a compact set $F\subset E_{\mathrm{unrect}}$ and Lipschitz
graph~$\Gamma$ such that ${\mathcal{H}^1(\Gamma\cap F)}\ge C\cdot
\mathcal{H}^1(F)$ that contradicts the fact that $F$ is purely
unrectifiable.

Let us now say a few words about the proof of Theorem~\ref{CMPT_theorem} given in~\cite{CMPT}.
Recall that this theorem is an analogue of Theorem~\ref{theorem_Leger}, where
the kernel $1/z$ and curvature squared $c^2(\mathcal{H}^1\lfloor E)$
are replaced by the kernels
$\kappa_n(z)=(\Re z)^{2n-1}/|z|^{2n}$, $n\in \mathbb{N}$, and corresponding permutations $p_{\kappa_n}$ (see (\ref{kappa_n})--(\ref{p_K(mu)})).
We will use the definitions given near the formula (\ref{O_tau}) and in the discussion of Theorem~\ref{theorem_Leger} above.
First we mention that  it is proved in \cite{CMPT} that the permutations $p_{\kappa_n}(z_1,z_2,z_3)$ behave similarly to $c^2(z_1,z_2,z_3)$ for all triangles with comparable
sides, whose one side makes a big angle with the vertical line. More precisely (see Lemma~\ref{lemma_CMPT_t=0} above),
it is shown there that for a fixed $\alpha_0\in(0,\pi/2)$ and given  $(z_1,z_2,z_3)\in
\mathcal{O}_\tau$, if the condition $(\ref{delta1})$, i.e.
$$
\theta_{V}(L_{z_1,z_2})+\theta_{V}(L_{z_2,z_3})+\theta_{V}(L_{z_1,z_3})\ge
\alpha_0,
$$
or
$$
\theta_{H}(L_{z_1,z_2})+\theta_{H}(L_{z_2,z_3})+\theta_{H}(L_{z_1,z_3})\le
\tfrac{3}{2}\pi-\alpha_0,
$$
is satisfied, then
$$
 p_{\kappa_n}(z_1,z_2,z_3)\ge C(\alpha_0,\tau)\cdot c(z_1,z_2,z_3)^2,\qquad C(\alpha_0,\tau)>0.
$$

This enables the authors of \cite{CMPT} to use the above-described
scheme from \cite{L} in order to construct the required Lipschitz graph
$\Gamma$ in the case when the first approximating line $D_0$ for
$\textrm{spt}\, \mu$ is far from the vertical line. Note that
exchanging the curvature for the permutations $p_{\kappa_n}$ still
requires new arguments in several key points of the proof.
Otherwise, when $D_0$ is close to the vertical line and the scheme
from \cite{L} does not work (as $\mu(F_3)$ may be too big), they
tune thresholds and apply some coverings so that they can use the result
for $D_0$, being far from the vertical line, to construct countably many Lipschitz graphs, which give $\Gamma$ after
appropriate joining.

We are now at the position to finish the proof of our Theorem~\ref{new_theorem2}.
This will be an adaptation of the arguments from \cite{CMPT}.

On the one hand, by the clause $(i)$ of Lemma~\ref{lemma_ends}, for a fixed $\alpha_0\in(0,\pi/2)$ and given  $(z_1,z_2,z_3)\in
\mathcal{O}_\tau$, if the condition (\ref{delta1}), i.e. the same as in the result for $t=0$ from \cite{CMPT} mentioned above,
is satisfied and $t=2-\tfrac{N}{n}$ for $n<N\le 2n$ or $t=\rho_{n,N}$ for $N\ge 2n$, then we also have
\begin{equation}
\label{last_section_ineq}
p_{K_t}(z_1,z_2,z_3)\ge C(\alpha_0,\tau)\cdot c(z_1,z_2,z_3)^2,\qquad C(\alpha_0,\tau)>0.
\end{equation}
It means that we can undeviatingly follow the scheme from
\cite{CMPT} (exchanging $p_{\kappa_n}$ for $p_{K_t}$) in order to
get our result for $t=2-\tfrac{N}{n}$, $n<N\le 2n$, and
$t=\rho_{n,N}$, $N\ge 2n$.

On the other hand, by the clause $(ii)$ of Lemma~\ref{lemma_ends},
the inequality~(\ref{last_section_ineq}) is true for
\begin{equation}
\label{last_t}
\begin{array}{ll}
t=-\tfrac{1}{2}\left(3+\sqrt{9-4\tfrac{N}{n}}\right), & n<N\le 2n, \\
t=-\frac{1}{2}\left(3+\rho_{n,N}+\sqrt{(3+\rho_{n,N})^2-4\frac{N}{n}}\right), & N\ge 2n,
\end{array}
\end{equation}
  only if both the conditions $(\ref{delta1})$ and $(\ref{delta2})$ are satisfied, i.e.
$$
\alpha_0\le \theta_{V}(L_{z_1,z_2})+\theta_{V}(L_{z_2,z_3})+\theta_{V}(L_{z_1,z_3})\le
\tfrac{3}{2}\pi-\alpha_0,
$$
or
$$
\alpha_0\le \theta_{H}(L_{z_1,z_2})+\theta_{H}(L_{z_2,z_3})+\theta_{H}(L_{z_1,z_3})\le
\tfrac{3}{2}\pi-\alpha_0,
$$
and thus the triangles $(z_1,z_2,z_3)$ are far from both the
vertical and horizontal line.

Consequently,  the scheme from \cite{CMPT} cannot be applied
directly for $t$ from (\ref{last_t}). However, as we will see, it
works after a few modifications (besides the exchange of
$p_{\kappa_n}$ for $p_{K_t}$) connected basically with adapting
geometrical arguments to both the conditions $(\ref{delta1})$ and
$(\ref{delta2})$. Since the cases where we are close to either the
vertical or horizontal line are well-separated and similar
geometrically, the arguments for the first approximating line $D_0$,
being  close (far) to (from) the vertical line, can be easily
transferred into the ones for $D_0$, being close (far) to (from) the
horizontal line.

We now reproduce the main steps of the proof, stemming from
\cite{CMPT}, with necessary changes when our permutations and the
conditions (\ref{delta1}) and (\ref{delta2}) are involved. Several
statements are formulated without proofs because they are the same
as in \cite{CMPT} (or \cite{L}) modulus the notation of
permutations.

Below we consider only $t$ from (\ref{last_t}). The
following two propositions will then imply
Theorem~\ref{new_theorem2} by the same contradiction arguments as in the proof
from \cite{L} (see the arguments around (\ref{contradiction}) above).
 Note that one has to take $\mu=40\mathcal{H}^1\lfloor F$ in
 Proposition~\ref{proposition2}, where the set $F$ is from Proposition~\ref{proposition1} (it may be suitably rescaled if necessary).

\begin{proposition}[An analogue of Lemma 3.4 in \cite{CMPT} and Proposition 1.1 in \cite{L}]
\label{proposition1}
Let $E\subset \mathbb{C}$ be a Borel set with $0<\mathcal{H}^1(E)<\infty$ and $p_{K_t}(\mathcal{H}^1\lfloor E)<\infty$.
Then for all $\eta>0$ there exists a set $F\subset E$ such that
\begin{itemize}
  \item $F$ is compact,
  \item $p_{K_t}(\mathcal{H}^1\lfloor F)\le \eta \diam F$,
  \item $\mathcal{H}^1(F)>\tfrac{1}{40}\diam F$,
  \item for all $z\in F$, for all $r>0$, $\mathcal{H}^1(F\cap B(z,r))\le 3r$.
\end{itemize}
\end{proposition}
\begin{proposition}
\label{proposition2}
For any constant $C_0\ge 10$, there exists a number $\eta>0$ such that if $\mu$ is any positive Radon measure on $\mathbb{C}$ satisfying
\begin{itemize}
  \item $\mu(B(0,1))\ge 1$, $\mu(\mathbb{C}\setminus B(0,2))=0$,
  \item for any ball $B$, $\mu(B)\le C_0 \diam B$,
  \item $p_{K_t}(\mu)\le \eta$,
\end{itemize}
then there exists a Lipschitz graph $\Gamma$ such that
$\mu(\Gamma)\ge 10^{-5}\mu(\mathbb{C})$.
\end{proposition}

The rest of the paper is devoted to the proof of Proposition~\ref{proposition2}, which is an analogue of Proposition~3.1 from \cite{CMPT}. First, we give several definitions (we confine ourselves to those which will be needed below; see \cite{CMPT,L} for further ones).
Let $\mu$ be a positive Radon measure on $\mathbb{C}$ and
$$
p_{K_t,\tau}(\mu)=\iiint_{\mathcal{O}_\tau}p_{K_t}(z_1,z_2,z_3)d\mu(z_1)d\mu(z_2)d\mu(z_3),
$$
see (\ref{p_K(mu)}) and (\ref{O_tau}).
For a ball $B=B(x,r)$ set
$$
\delta_\mu(x,r)=\frac{\mu(B(x,r))}{r}.
$$
We will use a small density threshold $\delta>0$ for this quantity.

 Given a fixed constant $k>1$, for any ball $B=B(x,r)\subset \mathbb{C}$
and $D$ a line in $\mathbb{C}$, set
\begin{align*}
&\beta_{1,\mu}^D(x,r)=\frac{1}{r}\int_{B(x,kr)}\frac{\dist(y,D)}{r}\;d\mu(y),\\
&\beta_{2,\mu}^D(x,r)=\left(\frac{1}{r}\int_{B(x,kr)}\left(\frac{\dist(y,D)}{r}\right)^2\;d\mu(y)\right)^{1/2}.
\end{align*}

Geometrical notation connected with lines and angles is given near the formula~(\ref{O_tau}) at the end of Section~\ref{section-auxiliary}.

\begin{lemma}
\label{Theorem_proof_lemma_1}
Let $\mu$ be a measure with linear growth $($with a constant $C_0$$)$, and $B(x,r)\subset
\mathbb{C}$ a ball with $\delta_\mu(x,r)\ge \delta$. Suppose that
$\tau$ is big enough, then for any $\varepsilon>0$, there exists
some $\delta_1=\delta_1(\delta,\varepsilon)>0$ such that
$$
\frac{p_{K_t,\tau}(\mu\lfloor kB)}{\mu(B)}\le \delta_1 \quad
\Longrightarrow \quad \inf_D\beta_{2,\mu}^D(x,r)\le\varepsilon.
$$
\end{lemma}
\begin{Proof}
The proof is the same as for Lemma 4.4 in \cite{CMPT}. We just have to use our Lemma~\ref{lemma_ends} instead of Lemma~2.3 there for the case when both the conditions (\ref{delta2}) and (\ref{delta1}) are satisfied,
and say that in the case
$$
\theta_H(L_{z_1,z_2})+\theta_H(L_{z_1,w})+\theta_H(L_{z_2,w})\le \alpha_0
 $$
we obtain the same estimate for $\textrm{dist}(w,L_{z_1,z_2})$ as in the case
 $$
\theta_V(L_{z_1,z_2})+\theta_V(L_{z_1,w})+\theta_V(L_{z_2,w})\le \alpha_0.
 $$
\end{Proof}

By Lemma~\ref{Theorem_proof_lemma_1}, chosen a point $x_0\in
\textrm{spt}\, \mu$, there exists an approximating line $D_0$ such
that $\beta_{1,\mu}^{D_0}(x_0,1)\le \varepsilon$. The next step is
to construct a first Lipschitz graph in the case when $D_0$ is far
from both the horizontal and vertical lines.

To do so, one first has to introduce a family of stopping time
regions and obtain the partition $\textrm{spt}\, \mu=Z\cup F_1\cup
F_2\cup F_3$ (see the exposition of the proof from \cite{L} above).
As this entirely repeats the corresponding part of \cite[Section
5]{CMPT} (cf. \cite[Subsection 3.1]{L}), we omit it. We just have to
mention that the thresholds $\theta_0$ and $\alpha$, arising there,
have to be adapted to that $D_0$ is far from both the horizontal and
vertical line. Namely, $\theta_0$ is now a threshold for both
$\theta_V(D_0)$ and $\theta_H(D_0)$. It means that one has to
distinguish not only the cases $\theta_V(D_0)\ge \theta_0$ and
$\theta_V(D_0)<\theta_0$ but also $\theta_H(D_0)\ge \theta_0$ and
$\theta_H(D_0)<\theta_0$. Moreover, $\alpha$ is tuned as follows: if
$\theta_V(D_0)$ or $\theta_H(D_0)$ are greater than $\theta_0$, then
$\alpha\le\theta_0/10$; if $\theta_V(D_0)$ or $\theta_H(D_0)$ are
not greater than $\theta_0$, then $\alpha=10\theta_0$.

Furthermore, see \cite{L,CMPT} for the way how one can define the Lipschitz function $A$ on the line $D_0$, using  $Z$, $F_1$, $F_2$, $F_3$, and appropriate thresholds.

Now we come to the main step of the proof of Proposition~\ref{proposition2}. The following lemma
is an analogue of Lemma~6.1 from \cite{CMPT}.
\begin{lemma}
\label{Theorem_proof_lemma_2}
Under the assumptions of Proposition~\ref{proposition2}, if furthermore
$$
\theta_0<\theta_V(D_0)<\tfrac{\pi}{2}-\theta_0,
$$
then there exists a Lipschitz graph $\Gamma$ such that $\mu(\Gamma)\ge \tfrac{99}{100}\mu(\mathbb{C})$.
\end{lemma}
For the proof one uses the above-mentioned function $A$ to obtain
the graph $\Gamma$, $Z\subset \Gamma$, and show that
$$
\mu(F_1)+\mu(F_2)+\mu(F_3)\le \tfrac{1}{100}\mu(\mathbb{C}).
$$
Indeed, the following lemmas are valid (recall that
$\mu(\mathbb{C})\ge 1$ by the assumptions).
\begin{lemma}
\label{Theorem_proof_lemma_4}
Under the assumptions of Proposition~\ref{proposition2},
$$
\mu(F_1)\le 10^{-6}.
$$
\end{lemma}
\begin{Proof}
This is an analogue of \cite[Proposition 6.3]{CMPT}, whose proof includes consideration of the two cases: 1) $\theta_V(D_0)>\theta_0$ and 2) $\theta_V(D_0)\le \theta_0$ (see the proof of \cite[Lemma~6.4]{CMPT}).

Under our settings, we have to consider three cases. Namely, the
case 1) has to be exchanged for
$\theta_0<\theta_V(D_0)<\tfrac{\pi}{2}-\theta_0$, although the proof
remains the same. The case 2) splits up into the following two:
$\theta_V(D_0)\le \theta_0$ and $\theta_V(D_0)\ge
\tfrac{\pi}{2}-\theta_0$ (i.e. $\theta_H(D_0)\le \theta_0$).
Arguments in the latter case are the same as in the former one.
\end{Proof}
\begin{lemma}[An analogue of Proposition 6.2 in \cite{CMPT}]
\label{Theorem_proof_lemma_3}
Under the assumptions of Proposition~\ref{proposition2},
$$
\mu(F_2)\le 10^{-6}.
$$
\end{lemma}
\begin{lemma}
\label{Theorem_proof_lemma_5}
Under the assumptions of Lemma~\ref{Theorem_proof_lemma_2},
$$
\mu(F_3)\le 10^{-6}.
$$
\end{lemma}
\begin{Proof}
The proof stems from the one of \cite[Proposition 6.5]{CMPT}, but
with exchange of $\theta_V(D_0)>\theta_0$ for
$\theta_0<\theta_V(D_0)<\tfrac{\pi}{2}-\theta_0$ as in
Lemma~\ref{Theorem_proof_lemma_2}.
\end{Proof}

\bigskip

Thus Proposition~\ref{proposition2} is proved under the assumptions
of Lemma~\ref{Theorem_proof_lemma_2}. What is left is to consider
the other case. The following statement satisfies the question.
\begin{lemma}
\label{Theorem_proof_lemma_6}
Under the assumptions of Proposition~\ref{proposition2}, if furthermore
$$
\theta_V(D_0)\le\theta_0\qquad \text{or} \qquad
\theta_V(D_0)\ge\tfrac{\pi}{2}-\theta_0 \quad (\text{i.e.}\quad
\theta_H(D_0)\le\theta_0),
$$
then there exists a Lipschitz graph $\Gamma$ such that
$\mu(\Gamma)\ge 10^{-5}\mu(\mathbb{C})$.
\end{lemma}
\begin{Proof}
To prove this, we repeat arguments from the proof of \cite[Lemma~7.1]{CMPT}, given for $\theta_V(D_0)\le\theta_0$, for the case $\theta_H(D_0)\le\theta_0$.
\end{Proof}

\section{Concluding remarks}

In this section we generalize Theorem~\ref{new_theorem2} to higher dimensions. Let us introduce necessary notation first.
For $d\in \mathbb{N}^+$ and $E\subset \mathbb{R}^d$ with finite length we consider the singular integral operator
$\mathbf{T}_{K_t} = (T^j_{K_t})^d_{j=1}$ such that formally
$$
T^j_{K_t}f(x):=\int_E f(y) K_t^j(x-y)d\mathcal{H}^1(y),\qquad K_t^j(x):=\kappa_N^j(x)+t\cdot \kappa_n^j(x),
$$
where
$\kappa_n^j(x):=x_j^{2n-1}/|x|^{2n}$ and $x=(x_1,\ldots,x_d)\in \mathbb{R}^d\setminus\{0\}$.
As before, we suppose that $N>n$, where $n,N\in \mathbb{N^+}$, and $t\in \mathbb{R}$. We also need the permutations
$$
\mathbf{p}_{K_t}(x,y,z):=\sum_{j=1}^d p_{K^j_t}(x,y,z) \qquad \text{ for distinct points }x,y,z\in \mathbb{R}^d,
$$
where $p_{K^j_t}(x,y,z)$ are the same as in (\ref{permutation_Main}) with $K^j_t$ instead of $K$. We also define the permutation of measure $\mathbf{p}_{K_t}(\mu)$ analogously to (\ref{p_K(mu)}).
\begin{theorem}
\label{new_theorem3}
Let $t$ be as mentioned in Theorem~\ref{new_theorem1}.
Given a Borel set $E\subset \mathbb{R}^d$ such that $0<\mathcal{H}^1(E)<\infty$,
if $\mathbf{p}_{K_t}(\mathcal{H}^1\lfloor E)<\infty$, then $E$ is rectifiable.
Moreover, if the operator $\mathbf{T}_{K_t}$ is $L^2(\mathcal{H}^1\lfloor E)$-bounded, then $E$ is rectifiable.
\end{theorem}

This result for $t=0$ was recently proved in \cite[Theorem 1.2(1) and Theorem 6.2]{Chousionis-Prat}.
To prove Theorem~\ref{new_theorem3} for all required $t$ we only need to use our Lemmas~\ref{lemma_new1},~\ref{lemma_new2} and \ref{lemma_ends} in order to show that for all $x,y,z\in \mathbb{R}^d$ such that $(x,y,z)\in \mathcal{O}_\tau$ and the assumptions of Lemma~\ref{lemma_ends} are satisfied,
\begin{equation}
\label{ineq_vectorial}
p_{K_t^j}(x,y,z)\ge C(t,\alpha_0,\tau)\, p_{K_0^j}(x,y,z),\qquad C(t,\alpha_0,\tau)>0.
\end{equation}
See the definitions of $\alpha_0$, $\tau$ and $\mathcal{O}_\tau$ before Lemma~\ref{lemma_ends1}. Then by  \cite[Proposition 3.3]{Chousionis-Prat}, adapted to the conditions (\ref{delta1}) and (\ref{delta2}), and the arguments similar to those in \cite[Section 6]{Chousionis-Prat} and our Section~\ref{section_proof} we immediately get the result. Note that \cite[Proposition 3.3]{Chousionis-Prat} slightly simplifies the approach from \cite{CMPT} (and improves Lemma~\ref{lemma_CMPT_t=0}) in the case $t=0$ as the parameter $\tau$ is not needed anymore. In our case this parameter is still necessary because of the inequality (\ref{ineq_vectorial}).

\bigskip

To finish, it is also worth mentioning here that under Ahlfors-David regularity
assumption one can expect that for $t$ as in Theorem~\ref{new_theorem1} the $L^2$-boundedness of the operator
associated with $K_t$ implies uniform rectifiability (the same for the $\mathbb{R}^d$ case). This is indeed
true. This result, among others with Ahlfors-David regularity
condition, has appeared in \cite{CMT}.

\section*{Acknowledgments}

I would like to express my sincere gratitude to Joan Mateu and
Xavier Tolsa for suggesting the problem and for many stimulating
conversations. I am also grateful to the Referee for
his/her valuable recommendations.

\end{document}